\documentclass[12pt, reqno]{amsart}
\usepackage[foot]{amsaddr}

\usepackage[mathscr]{eucal}
\usepackage{graphicx}
\usepackage{caption}
\usepackage{subcaption}
\usepackage{multicol}
\usepackage{amssymb,amsfonts,amsthm}
\usepackage{color}
\usepackage[centertags]{amsmath}
\usepackage{empheq}
\usepackage{geometry}
\usepackage{algorithm}
\usepackage{placeins}
\usepackage{algpseudocode}
\usepackage{hyperref}
\usepackage{mathrsfs}
\usepackage{enumerate}
\usepackage{esint}
\usepackage{booktabs}
\usepackage{multirow}
\usepackage{cite}
\usepackage{bm}
\usepackage{tikz}
\usetikzlibrary{arrows.meta,positioning,calc,fit,backgrounds}
\definecolor{purple}{rgb}{0.5,0,0.5}

\theoremstyle{plain}
\newtheorem{thm}{Theorem}[section]
\newtheorem{prop}[thm]{Proposition}

\newtheorem{remark}[thm]{Remark}

\numberwithin{equation}{section}
\numberwithin{figure}{section}

\newcommand{\R}{\mathbb{R}}
\newcommand{\N}{\mathbb{N}}
\newcommand{\E}{\mathbb{E}}
\newcommand{\Pcal}{\mathcal{P}}

\newcommand{\Lcal}{\mathcal{L}}
\newcommand{\Ecal}{\mathcal{E}}
\newcommand{\Fcal}{\mathcal{F}}
\newcommand{\Hmone}{\dot{H}^{-1}(\Omega)}
\newcommand{\bv}{\mathbf{v}}
\newcommand{\bn}{\mathbf{n}}
\newcommand{\be}{\mathbf{e}}

\newcommand{\bX}{\mathbf{X}}

\newcommand{\blam}{\bm{\lambda}}

\title[Self-Supervised In-Context Operator Learning for Optimal Mixing]{Self-Supervised In-Context Operator Learning for Optimal Mixing in Incompressible Flows}

\author{Weiwei Hu$^1$}
\address{$^1$Department of Mathematics, University of Georgia, Athens, GA 30602, USA}
\author{Rongjie Lai$^2$}
\address{$^2$Department of Mathematics, Purdue University, 150 N.~University Street, West Lafayette, IN 47907, USA}
\author{Ziqian Li$^3$}
\address{$^3$Chair for Dynamics, Control, Machine Learning and Numerics, Department of Mathematics, Friedrich-Alexander-Universit\"{a}t Erlangen-N\"{u}rnberg, Erlangen 91058, Germany}

\email{Weiwei.Hu@uga.edu, lairj@purdue.edu, ziqian.li@fau.de}

\begin{document}

\begin{abstract}
We develop a self-supervised operator-learning framework for optimal mixing in incompressible flows, aimed at overcoming the computational limitations of classical PDE-constrained optimal control. The objective is to design velocity fields that efficiently homogenize a passive scalar while minimizing kinetic energy, quantified through the $\dot H^{-1}(\Omega)$ mix-norm. 
Classical methods solve a separate optimal control problem for each initial condition using iterative forward–adjoint procedures, leading to substantial computational costs and limited scalability.
In contrast, we propose to learn the solution operator that maps initial scalar fields to optimal mixing controls directly. The proposed self-supervised method does not require access to supervised labels, thereby significantly reducing the data-generation burden of supervised operator-learning methods.  Our framework integrates a Lagrangian reformulation of the transport equation with a transformer-based neural operator architecture. The transport dynamics are represented through characteristic ODEs, enabling particle-based Monte Carlo optimization, while the learned velocity field is parameterized in a divergence-free spectral basis enforcing incompressibility and no-penetration conditions exactly. Numerical experiments on parametric families of initial conditions in two and three spatial dimensions demonstrate robust generalization to unseen regimes and approximately exponential decay of the mix-norm.  The learned controls also transfer consistently to independent Eulerian finite-volume simulations, validating the physical fidelity of the learned operator. To the best of our knowledge, this work provides the first self-supervised operator-learning framework for optimal mixing and among the first data-driven approaches for optimal mixing in fully three-dimensional incompressible flows.

\end{abstract}

\maketitle

%% ====================================================================
\section{Introduction}\label{sec:intro}

Optimal mixing of passive scalars by incompressible flows is a fundamental problem at the intersection of fluid dynamics, control theory, and scientific computing, with applications ranging from microfluidics and chemical processing to geophysical transport and atmospheric dynamics \cite{stroock2002chaotic}. Fluid mixing is governed by the interplay between advection and diffusion \cite{villermaux2019mixing,thiffeault2012using}. In many applications characterized by large Péclet numbers, diffusive effects are negligible over the time scales of interest, and the mixing process is driven primarily by advection \cite{lin2011optimal,alberti2016exponential}. In this regime, stirring-induced stretching and folding become the dominant mechanisms responsible for generating fine-scale structures and enhancing mixing.

 In the absence of molecular diffusivity, the mathematical formulation of optimal mixing typically involves the transport equation coupled with incompressibility constraints and seeks velocity fields that enhance mixing under constraints on energy, enstrophy, or action (e.g., \cite{balasuriya2015dynamical, cortelezzi2008feasibility, couchman2010control, crippa2017cellular, eggl2020mixing, eggl2022mixing, franjione1992symmetry, gubanov2010towards, lin2011optimal, liu1994quantification, lunasin2012optimal, mathew2007optimal, sharma1997control, thiffeault2012using, vikhansky2002control, vikhansky2002enhancement}). A common objective is to optimize a quantitative measure of mixing at a prescribed terminal time while controlling the cost of stirring.

% We now introduce the mathematical framework underlying this class of optimal mixing problems. Let $\Omega\subset\mathbb{R}^d$, $d=2,3$, be a bounded domain with sufficiently regular boundary $\partial{\Omega}$, and let 
%\[
%\bv:[0,T]\times \Omega \to \mathbb{R}^d
%\]
%be an incompressible velocity field satisfying
%\begin{align}
%\nabla \cdot \bv =0
%\quad \text{in } \Omega,
%\quad \text{and}\quad
%\bv \cdot \bn =0
%\quad \text{on } \partial{\Omega}, \label{prop_v}
%\end{align}
%where $\bn$ denotes the outward unit normal vector on $\partial{\Omega}$. Given an initial scalar distribution $\rho_0$, the transported scalar field $\rho$ is governed by
%\begin{equation}
%\partial_t \rho + \bv \cdot \nabla \rho =0,
%\qquad
%\rho(0)=\rho_0.
%\label{eq:transport}
%\end{equation}
%Using the incompressibility  and no-penetration boundary  condition in \eqref{prop_v}, it is easy to verify that any $L^p$-norm of $\rho$ is conserved, i.e.,
%\begin{align}
%\|\rho(t)\|_{L^p(\Omega)}=\|\rho_{0}\|_{L^p(\Omega)}, \quad t\geq 0, \quad p\in[1,\infty]. \label{theta_Lp}
%\end{align}

A central question in optimal  mixing is how to design velocity fields that efficiently redistribute scalar variance from low to high frequencies, thereby generating progressively finer filaments. More specifically, 
we consider the optimal mixing of two fluid species in a domain
$\Omega \subset\mathbb{R}^d$, $d=2,3$ with sufficiently regular boundary $\partial{\Omega}$. The species are transported by an incompressible velocity field
$\mathbf{v}(t,\mathbf{x})$ satisfying:
\begin{align}
\nabla \cdot \bv =0
\quad \text{in } \Omega,
\quad \text{and}\quad
\bv \cdot \bn =0
\quad \text{on } \partial{\Omega}, \label{prop_v}
\end{align}
where $\bn$ denotes the outward unit normal vector on $\partial{\Omega}$.
To distinguish the two species, we consider a \emph{signed color field}
$\rho(t,\mathbf{x})$, taking values $+1$ and $-1$ for the two species,
respectively. Unlike a probability density, the signed color field is
not required to be nonnegative. Instead, it characterizes the spatial
distribution of the two species.
The color field evolves according to the transport equation
\begin{equation}
    \partial_t \rho + \mathbf{v}\cdot\nabla\rho = 0,
    \qquad \nabla\cdot\mathbf{v}=0,
    \label{eq:transport}
\end{equation}
with initial condition $\rho(0,\mathbf{x})=\rho_0(\mathbf{x})$.

Following the classical works of Mathew, Mezi\'c, and Petzold in \cite{mathew2005multiscale} and subsequent developments in the mathematical mixing literature (e.g.~\cite{alberti2016exponential}), we use the homogeneous Sobolev norm $\dot H^{-1}(\Omega)$ as a mix-norm. Smaller values of $\|\rho\|_{\dot H^{-1}(\Omega)}$ correspond to better mixing through the suppression of large-scale structures.
Our objective is to determine an incompressible velocity field that
promotes mixing while minimizing the associated control cost.
More specifically, given an initial datum $\rho_0\in L^\infty(\Omega)$, the \emph{optimal mixing problem} takes the form
\begin{equation*}
\min_{\bv} \mathcal{J}(\mathbf{v}) :=\ \lambda_{\mathrm{kin}}  \int_0^T \int_\Omega
    \|\mathbf{v}(t,\mathbf{x})\|^2
    \,d\mathbf{x}\,dt + \lambda_{\mathrm{mix}}\|\rho(T,\cdot)\|_{\dot H^{-1}(\Omega)}^{2}
\end{equation*}
subject to \eqref{prop_v} and \eqref{eq:transport},
where $T>0$ is a prescribed terminal time. The first term penalizes the kinetic energy of the stirring velocity,
while the second measures the degree of mixing through the negative
Sobolev norm, which quantifies large-scale spatial inhomogeneity and
decreases as the two species become increasingly mixed.
Importantly, the kinetic energy is evaluated with respect to the
spatial volume measure rather than the signed color field, ensuring
that the control cost remains nonnegative.

Classical approaches formulate optimal mixing as a PDE-constrained optimal control problem in the Eulerian frame. One widely used formulation minimizes a weighted combination of the terminal mix-norm and a regularization term on the velocity field, subject to the transport equation \eqref{eq:transport} and conditions in \eqref{prop_v} (e.g.~\cite{mathew2007optimal,hu2025cellular, hu2026structure}). Such formulations have led to important analytical and computational developments, including variational formulations, adjoint-based gradient methods, optimal stirring protocols, and structure-preserving numerical schemes. However, these approaches suffer from a major computational bottleneck: each realization of the initial condition $\rho_0$ requires solving a new forward--adjoint optimization problem, often on highly refined meshes due to the multiscale nature of mixing (e.g.~\cite{zheng2023numerical}).

This limitation becomes particularly severe in many-query settings, including uncertainty quantification \cite{smith2013uncertainty}, real-time flow control \cite{garcia1989model}, inverse design \cite{borrvall2003topology}, and optimization over large families of initial data \cite{benner2015survey}. In such scenarios, repeatedly solving PDE-constrained optimization problems quickly becomes computationally prohibitive, especially in three-dimensional settings. Recent works have sought to alleviate these difficulties through  feedback control designs  and structure-preserving discretization, etc. \cite{hu2023feedback,hu2026structure}. However, existing methods remain essentially \emph{instance-based}: the optimal control must still be recomputed whenever the initial scalar distribution changes.

\subsection{A New Paradigm: Learning the Control Operator}

Motivated by recent work on self-supervised in-context learning for mean-field games \cite{huang2025unsupervised}, we introduce a fundamentally different computational paradigm for optimal mixing. Rather than repeatedly solving an optimization problem for each initial condition, we directly \emph{learn the solution operator}
\[
\rho_0 \mapsto \bv,
\]
mapping initial scalar fields to approximately optimal mixing controls. This viewpoint shifts the problem from repeated PDE-constrained optimization to operator learning. Once trained offline, the learned operator produces a mixing control for a new initial condition through a single forward evaluation, enabling near real-time inference.

Our framework combines two main ingredients. First, we reformulate the transport dynamics in Lagrangian coordinates. Instead of working directly with the Eulerian transport equation, we describe the dynamics using characteristic trajectories, which satisfy a system of ordinary differential equations for particle motion. This reformulation enables Monte Carlo approximation of the objective functional and allows gradients to be computed through explicit time integration without requiring adjoint PDE solvers.
Second, we develop a self-supervised, transformer-based in-context neural operator that encodes the initial scalar distribution as point clouds corresponding to the two phases of the mixture. Unlike supervised operator-learning methods that rely on solver-generated input--output pairs, our method learns directly from the transport dynamics and control objective, without requiring reference controls computed by another numerical solver. The network outputs time-dependent coefficients of a divergence-free spectral basis, guaranteeing incompressibility and no-penetration boundary conditions exactly by construction. The resulting architecture combines physical structure preservation with the expressive power of modern operator-learning methods. To verify the learned solution operator, we use an independent finite-volume solver to evolve the Eulerian transport equation under the learned velocity field. The scalar-field patterns and fitted mix-norm decay rate from the Eulerian simulation agree with the corresponding Lagrangian particle-based results (Section~\ref{sec:verification}).

We summarize main contributions of this work as follows.

\begin{enumerate}[(1)]
\item \textbf{A Lagrangian reformulation of optimal mixing.} We represent the two species by particles carrying fixed labels $+1$ or $-1$, assigned according to their initial locations and transported along the flow. The flow map generated by $\bv$ expresses the kinetic energy and terminal mix-norm in terms of particle trajectories, with the signed labels entering the mix-norm. Replacing the transport PDE with characteristic ODEs allows us to estimate the cost by Monte Carlo and compute its gradients by differentiating through an explicit time integrator, without solving a separate adjoint PDE (Section~\ref{sec:lagrangian}).

\item \textbf{A self-supervised operator-learning formulation.} We learn the control operator $\rho_0\mapsto\bv^{\ast}$, which maps an initial scalar field to an approximately optimal stirring control. The operator is trained directly using the transport dynamics and control objective, without solver-generated optimal-control labels. Once trained, it produces a control for a new initial condition with a single forward evaluation (Section~\ref{sec:operator}).

\item  \textbf{A physics-constrained transformer architecture.} We parameterize the velocity in a divergence-free spectral basis that enforces incompressibility and no-penetration exactly. A transformer with bidirectional cross-attention encodes the two phases of the initial density as point clouds and outputs time-dependent basis coefficients. We train the operator by backpropagating through the ODE rollout (Sections~\ref{sec:divfree_basis}--\ref{sec:context_encoder} and~\ref{sec:training}).

\item \textbf{Generalization across parametric families in 2D and 3D.} We evaluate four parametric families of two-phase densities in 2D and four in 3D, each with a held-out region of parameter space. On held-out data, the learned operator maintains its terminal mix-norm performance and fitted decay rate. To our knowledge, this is the first self-supervised operator-learning framework for optimal mixing and one of the first data-driven approaches to optimal mixing in fully three-dimensional incompressible flows (Sections~\ref{sec:dataset} and~\ref{sec:experiments}).

%\item \textbf{Verification against the Eulerian transport PDE.} We use an independent finite-volume solver to simulate the Eulerian transport equation with the learned velocity field. The resulting scalar-field structure and fitted mix-norm decay rate agree with those obtained from the Lagrangian particle simulation (Section~\ref{sec:verification}).
\end{enumerate}

\subsection{Related works}

Our work connects and extends three major research directions.

\medskip

\noindent
\textbf{Mixing and optimal control.}
Optimal mixing and stirring enhancement have been studied extensively through variational formulations, optimal-transport perspectives, cellular flows, and constrained stirring strategies. Existing approaches include mix-norm-based optimal-control formulations, least-action principles, enstrophy-constrained stirring, and finite-dimensional flow design (e.g.~\cite{mathew2007optimal, lin2011optimal, hu2025cellular, hu2026structure, li2026hamiltonian}). These methods have yielded important analytical insights into mixing rates, filamentation mechanisms, and energy-constrained transport. However, they primarily focus on solving an individual optimization problem for each fixed initial condition.

\medskip

\noindent
\textbf{Feedback control and computation.}
To reduce computational cost, prior work has explored feedback-control laws, instantaneous control strategies, and structure-preserving numerical discretizations (e.g.~\cite{hu2023feedback,zheng2023numerical, hu2026structure}). Recent developments have emphasized preserving physical invariants, maintaining consistency between state and adjoint discretizations, and enforcing conservation properties at the discrete level (e.g.~\cite{hu2026structure}). Although these methods improve numerical robustness and computational efficiency, they remain tied to repeated PDE solves and do not address the many-query demands of optimal-mixing design.

\medskip

\noindent
\textbf{Operator learning and neural operators.}
Neural operators \cite{kovachki2023neural} approximate mappings between function spaces. Prominent architectures include DeepONet \cite{lu2021deeponet}, the Fourier neural operator \cite{li2020fourier}, and In-Context Operator Networks \cite{incontext2023}. Many such methods are trained in a supervised manner and therefore require input-output pairs, often generated by numerical solvers. In-Context Operator Networks represent input functions through finite sets of point samples, allowing them to accommodate varying sample locations and numbers. This set-based processing draws on permutation-invariant representations in the spirit of Deep Sets \cite{zaheer2017deepsets} and on attention mechanisms \cite{vaswani2017attention}. Self-supervised in-context operator learning has been developed for optimal transport and mean-field games~\cite{huang2025unsupervised}, second-order mean-field control \cite{gao2026self}, and Wasserstein gradient flows \cite{feng2026learn}, and has been analyzed theoretically for optimal transport on probability-measure spaces~\cite{cole2026context}. These self-supervised approaches do not require precomputed numerical solutions as training labels. Our contribution relative to these works is a two-density context encoder with bidirectional cross-attention, a divergence-free output basis that enforces the physical constraints exactly, and a Lagrangian particle loss tailored to the mix-norm objective.

Our work lies at the intersection of these areas and, to the best of our knowledge, provides the first self-supervised operator-learning framework for optimal mixing control in incompressible flows.

\subsection{Organization of the Paper}

The remainder of the paper is organized as follows. In Section~\ref{sec:lagrangian}, we introduce the $\dot H^{-1}$ mix-norm and derive the Lagrangian reformulation of the optimal mixing problem. Section~\ref{sec:operator} presents the neural operator framework, including the divergence-free spectral basis and transformer-based architecture. Section~\ref{sec:training} describes the particle-based training methodology and optimization procedure, and Section~\ref{sec:dataset} introduces the parametric families of initial conditions used in the numerical experiments. Section~\ref{sec:experiments} reports two- and three-dimensional numerical results together with Eulerian validation studies. Finally, Section~\ref{sec:conclusion} concludes the paper and discusses future research directions.

%% ====================================================================
\section{Mix-norm and Lagrangian formulation}\label{sec:lagrangian}

This section defines the $\Hmone$ mix-norm through its dual, elliptic, and Green kernel representations. We then use the volume-preserving flow map to express the control objective in Lagrangian coordinates and derive its particle approximation. The truncated Fourier expansion used to evaluate the mix-norm is given in Section~\ref{sec:loss_computation}.

\subsection{The $\dot H^{-1}$ mix-norm}\label{sec:mixnorm}
We represent the two fluid species in $\Omega\subset\mathbb{R}^d$ by a signed color field. Let $A\subset\Omega$ be measurable with $0<|A|<|\Omega|$, and write $A^c=\Omega\setminus A$. The initial field is
\begin{equation*}
\rho_0=\mathbf{1}_{A}-\mathbf{1}_{A^{c}},
\end{equation*}
where $+1$ and $-1$ label the first and second species, respectively. We sample the two phases from the uniform probability measures
\begin{equation*}
\rho_1=\mathrm{Uniform}(A)\qquad\text{and}\qquad\rho_2=\mathrm{Uniform}(A^c).
\end{equation*}
The field $\rho_0$ records the initial spatial arrangement and is not a probability density. Incompressible advection without diffusion preserves the values $+1$ and $-1$ of $\rho(t,\cdot)$. The two phases need not have equal volume, so the spatial mean need not vanish.

The multiscale mix-norm introduced by Mathew, Mezi\'c, and Petzold \cite{mathew2005multiscale} motivates this choice. Here, following the optimal-stirring formulation of Lin et al. \cite{lin2011optimal}, we use the homogeneous Sobolev seminorm $\Hmone$. For mean-zero $\rho\in L^{2}_0(\Omega):=\{\rho\in L^{2}(\Omega): \int_{\Omega}\rho\,dx=0\}$,  we define its homogeneous $\Hmone$ norm by \begin{equation}\label{eq:neg_norm}
\|\rho\|_{\Hmone} = \sup\Bigl\{\textstyle\int_\Omega \rho\phi\,dx: \phi\in H^{1}(\Omega),\ \int_\Omega\phi\,dx=0,\ \int_\Omega|\nabla\phi|^2\,dx\le 1\Bigr\}.
\end{equation} 
We then define \(\dot H^{-1}(\Omega)\) as the completion of \(L^2_0(\Omega)\) under this norm: \[
\dot H^{-1}(\Omega)
:=\overline{L^2_0(\Omega)}^{\,\|\cdot\|_{\dot H^{-1}}}.
\]
This defines a seminorm on $L^2(\Omega)$ and a norm on $L^2_0(\Omega)$. It measures large-scale inhomogeneity, i.e.  low-frequency components contribute more strongly than fine-scale components, so decay of the seminorm reflects redistribution of scalar content toward finer scales.

Our signed color field need not have zero mean, so we define its mix-norm by projecting out its spatial mean:
\begin{equation}\label{eq:P0_convention}
\|\rho\|_{\Hmone} := \|\mathcal{P}_{0}\rho\|_{\Hmone},
\qquad
\mathcal{P}_{0}\rho := \rho-\frac{1}{|\Omega|}\int_{\Omega}\rho\,dx,
\end{equation}
where $\mathcal{P}_{0}$ is the $L^{2}$-orthogonal projection onto $L^{2}_{0}(\Omega)$. This convention is consistent in time. Incompressibility and the no-penetration condition \eqref{prop_v} conserve the spatial mean $\bar{\rho}:=|\Omega|^{-1}\int_\Omega \rho(t,\mathbf{x})\,dx$. Thus, perfect mixing corresponds to weak convergence of $\rho(t,\cdot)$ to $\bar{\rho}$ together with $\|\mathcal{P}_{0}\rho(t,\cdot)\|_{\Hmone}\to0$.

The Neumann Green kernel below has zero mean in each argument and therefore removes the constant mode implicitly. The Fourier expansion in Section~\ref{sec:loss_computation} does the same by omitting $k=0$. The elliptic formulation, by contrast, must use $\mathcal{P}_{0}\rho$ on its right-hand side to satisfy the Neumann solvability condition.

For computation, we use two equivalent representations. First, let $\eta$ be the unique mean-zero solution of the Neumann problem
\begin{equation}
\label{eq:elliptic}
-\Delta \eta = \mathcal{P}_{0}\rho \text{ in } \Omega, \qquad \partial_{\bn}\eta = 0 \text{ on } \partial{\Omega}, \qquad \int_{\Omega}\eta\,dx = 0.
\end{equation}
When $\eta\ne0$, the supremum in \eqref{eq:neg_norm} is attained at $\phi=\eta/\|\nabla\eta\|_{L^2}$. If $\eta=0$, both sides of the identity below vanish. Therefore,
\begin{equation}\label{eq:Hminus1_identity}
\|\rho\|^2_{\Hmone} = \langle \rho,\eta\rangle_{L^{2}} = \int_\Omega |\nabla\eta|^{2}\,dx.
\end{equation}
Second, let $\mathcal{K}:\Omega\times\Omega\to\R$ be the Neumann Green kernel of $(-\Delta)^{-1}$ on $L^{2}_0(\Omega)$, characterized by
\begin{equation*}
-\Delta_x \mathcal{K}(x,y) = \delta_y(x) - \frac{1}{|\Omega|} \text{ in }\Omega, \qquad \partial_{\bn_x}\mathcal{K}(x,y) = 0 \text{ on }\partial{\Omega}, \qquad \int_\Omega \mathcal{K}(x,y)\,dx = 0.
\end{equation*}
Since $\eta(x) = \int_{\Omega}\mathcal{K}(x,y)\rho(y)\,dy$, substituting into \eqref{eq:Hminus1_identity} gives
\begin{equation}\label{eq:Hminus1_kernel}
\|\rho\|^2_{\Hmone} = \iint_{\Omega\times\Omega} \rho(x)\,\mathcal{K}(x,y)\,\rho(y)\,dx\,dy,
\end{equation}
which is the form we use under the Lagrangian change of variables.

For the computational framework developed in this paper, we specialize to the unit square or cube, \(\Omega=(0,1)^d\), with \(d\in\{2,3\}\). Section~\ref{sec:loss_computation} describes how we evaluate the kernel integral by a truncated Fourier expansion and particle sampling.

\subsection{From Eulerian to Lagrangian control}\label{sec:eul_to_lag}

The Eulerian formulation of optimal mixing penalizes the kinetic energy of $\bv$ over $[0,T]\times\Omega$ and the terminal mix-norm of $\rho$.

To incorporate incompressibility and the no-flux condition, we restrict the velocity field to a finite-dimensional admissible space and write
\begin{equation}\label{eq:vel_expansion}
\bv(t,x)=\sum_{\ell=1}^{L}\lambda_\ell(t)\,\be_\ell(x),
\end{equation}
where $\lambda_\ell:[0,T]\to\mathbb{R}$ are time-dependent control inputs. We assume the basis fields $\be_\ell$ are sufficiently regular and satisfy
\begin{align}\label{prop_basis}
\nabla\cdot\be_\ell=0 \quad\text{in }\Omega,
\qquad
\be_\ell\cdot\bn=0 \quad\text{on }\partial\Omega.
\end{align}

Let
$
V_L:=\operatorname{span}\{\be_1,\ldots,\be_L\}
$
denote the $L$-dimensional admissible velocity space introduced in \eqref{eq:vel_expansion}, where each $e_\ell\in W^{1,\infty}(\Omega), \ell=1,2,\dots,L,$ and satisfies \eqref{prop_basis}. Define
$$
\Fcal:=L^2(0,T;V_L)
=
\left\{
\bv(t,x)=\sum_{\ell=1}^L\lambda_\ell(t)\be_\ell(x):
\lambda_\ell\in L^2(0,T)
\right\}.
$$
Clearly,
$
\Fcal\subset L^1\bigl(0,T;W^{1,\infty}(\Omega;\mathbb R^d)\bigr).
$
Moreover, the divergence-free and no-flux conditions in \eqref{prop_v} are inherited from the basis fields. Thus, for every \(\bv\in\Fcal\), Carathéodory theory for ordinary differential equations yields a unique flow map \(\bX\), which is bi-Lipschitz in space for each \(t\in[0,T]\). This regularity provides the flow-map representation required for the Lagrangian reformulation.

Under this finite-dimensional parameterization, the Eulerian optimal-mixing problem is
\begin{equation}\label{eq:eulerian_opt}
\begin{aligned}
\min_{\bv\in\Fcal}\quad
&\lambda_{\mathrm{kin}}
\int_0^T\int_\Omega|\bv(t,x)|^2\,dx\,dt
+
\lambda_{\mathrm{mix}}
\iint_{\Omega\times\Omega}
\rho(T,x)\,\mathcal{K}(x,y)\,\rho(T,y)\,dx\,dy \\
&\text{subject to the transport equation \eqref{eq:transport}.}
\end{aligned}
\end{equation}
The flow map $t\mapsto\bX(t;x)$ is generated by $\bv$ through the characteristic ODE
\begin{equation}\label{eq:ode_traj}
\frac{d}{dt}\bX(t;x) = \bv(t,\bX(t;x)), \quad t\in(0,T]; \qquad \bX(0;x) = x.
\end{equation}
Along characteristics, $\rho(t,\bX(t;x))=\rho_0(x)$. Equivalently, the transported signed measure satisfies
\[
\rho(t,\cdot)\,dx = \bX(t,\cdot)_{\#}\bigl(\rho_0(x)\,dx\bigr).
\]
In particular, for any test function $h$,
\begin{equation}\label{eq:pushforward}
\int_\Omega h(x)\,\rho(t,x)\,dx = \int_\Omega h(\bX(t;x))\,\rho_0(x)\,dx.
\end{equation}

Both terms in the Eulerian objective can be expressed through the flow map. The kinetic-energy term does not involve the density. Since $\nabla\cdot\bv=0$ in $\Omega$ and $\bv\cdot\bn=0$ on $\partial\Omega$, the flow map preserves volume and maps $\Omega$ onto itself. In particular, its Jacobian $J(t,x)=\det D_x\bX(t;x)$ satisfies $\partial_tJ=(\nabla\cdot\bv)(t,\bX(t;x))J=0$ and $J(0,x)=1$, so $J(t,x)=1$ for almost every $x$ and every $t$. Thus, changing variables via $z=\bX(t;x)$ gives
\[
\int_0^T\int_\Omega |\bv(t,z)|^2\,dz\,dt
=
\int_\Omega
\int_0^T |\bv(t,\bX(t;x))|^2\,dt\,dx.
\]
Applying the pushforward identity \eqref{eq:pushforward} to both variables in the kernel representation \eqref{eq:Hminus1_kernel} yields
\[
\|\rho(T)\|_{\Hmone}^2
=
\iint_{\Omega\times\Omega}
\rho_0(x)\,
\mathcal{K}\bigl(\bX(T;x),\bX(T;y)\bigr)\,
\rho_0(y)\,dx\,dy.
\]
Consequently, the Eulerian problem \eqref{eq:eulerian_opt} is equivalent to the Lagrangian optimization problem
\begin{equation}\label{eq:lagrangian_opt}
\min_{\bv\in\Fcal}
\lambda_{\mathrm{kin}}\int_\Omega
\int_0^T |\bv(t,\bX(t;x))|^2\,dt\,dx
+
\lambda_{\mathrm{mix}}
\iint_{\Omega\times\Omega}
\rho_0(x)\,
\mathcal{K}\bigl(\bX(T;x),\bX(T;y)\bigr)\,
\rho_0(y)\,dx\,dy.
\end{equation}

\noindent\textbf{Particle discretization. } We approximate the integrals in \eqref{eq:lagrangian_opt} by uniform sampling on the unit square or cube. Since $|\Omega|=1$, the spatial volume measure is a probability measure. We draw independent and identically distributed samples
\[
 x_1,\ldots,x_N\overset{\mathrm{i.i.d.}}{\sim}\operatorname{Uniform}(\Omega),
\]
write $\bX_i(t):=\bX(t;x_i)$, and assign each particle its fixed color
\[
 w_i:=\rho_0(x_i)\in\{+1,-1\}.
\]
The kinetic term is an unweighted average over these particles. The mix-norm uses products of their signed weights, which we center to give the empirical signed measure zero total mass. We exclude the singular diagonal of $\mathcal K$ (Remark~\ref{rem:kernel_smoothing}) and obtain the particle objective
\begin{equation}\label{eq:lagrangian_is}
\min_{\bv\in\Fcal}\
\frac{\lambda_{\mathrm{kin}}}{N}\sum_{i=1}^N\int_0^T|\bv(t,\bX_i(t))|^2\,dt
\;+\;
\frac{\lambda_{\mathrm{mix}}}{N(N-1)}\sum_{\substack{i,j=1\\ i\ne j}}^N \tilde w_i\,\tilde w_j\,\mathcal{K}\bigl(\bX_i(T),\bX_j(T)\bigr),
\end{equation}
where $\tilde w_i:=w_i-\frac1N\sum_{j=1}^N w_j$ are the centered weights. They implement the projection $P_0$ of \eqref{eq:P0_convention} at the particle level: the signed empirical measure $\frac1N\sum_{i}\tilde w_i\,\delta_{\bX_i(T)}$ has total mass zero and is the particle analogue of $P_0\rho(T,\cdot)$. The off-diagonal sum in \eqref{eq:lagrangian_is} is a Monte Carlo approximation of the squared terminal mix-norm. Section~\ref{sec:training} replaces $\mathcal{K}$ by a smooth truncated kernel, for which the diagonal is finite and may be retained, and states the resulting training loss in closed form.

\section{Neural operator framework}\label{sec:operator}

We now approximate the control map $\rho_0\mapsto\bv^{\ast}$ by a neural operator trained on the Lagrangian objective. After defining the amortized loss, we specify its divergence-free output basis and a transformer architecture that encodes the two phases as point clouds and returns time-dependent coefficients.

\subsection{Solution operator and amortized objective}\label{sec:solution_op}
Solving \eqref{eq:lagrangian_opt} for a single $\rho_0$ produces an optimal velocity field $\bv^{\ast}$ tailored to that input. Operator learning aims to approximate the map $\rho_0\mapsto\bv^{\ast}$ by a single neural network. We represent each initial field by the pair of uniform probability measures $(\rho_1,\rho_2)$ defined in Section~\ref{sec:mixnorm}. Writing $\Pcal(\Omega)$ for the space of probability measures on $\Omega$, we therefore seek an operator
\begin{equation*}
G:\Pcal(\Omega)\times\Pcal(\Omega)\times[0,T]\longrightarrow\R^{L}, \qquad (\rho_1,\rho_2,t)\longmapsto \blam(t),
\end{equation*}
where $L$ is the number of elements in the truncated spectral basis (see Section~\ref{sec:divfree_basis}) and the output $\blam(t)=(\lambda_1(t),\ldots,\lambda_L(t))$ parameterizes the velocity through that basis $V_L$. The target space of $G$ is the finite-dimensional $\R^{L}$, independent of the grid, which makes amortization tractable. 

The two measures are presented to the network as two point clouds 
$$\{x_j^{(1)}\}_{j=1}^{n_1}\sim\rho_1\quad \text{and} \quad \{x_j^{(2)}\}_{j=1}^{n_2}\sim\rho_2,$$ with $n_1,n_2\in\N$,
 the number of samples drawn i.i.d.\ from each phase, and the pair is referred to as the context. The approximation $G_\Theta$ of $G$, parameterized by a vector $\Theta$ of network parameters, takes the context and a query time $t$ and returns the coefficient vector
\begin{equation*}
\blam(t) = G_\Theta\big(\{x_j^{(1)}\}_{j=1}^{n_1},\,\{x_j^{(2)}\}_{j=1}^{n_2};\,t\big).
\end{equation*}
Two structural properties are imposed by construction. First, $G_\Theta$ is permutation-invariant within each group, because the samples $\{x_j^{(k)}\}$ represent the measure $\rho_k$ and carry no order. Second, $G_\Theta$ is sampling-invariant: the output stabilizes as $n_1,n_2\to\infty$ to a functional of $\rho_1,\rho_2$. Both properties are enforced by mean-pooling over tokens, as described in Section~\ref{sec:context_encoder}.

To learn a single $G_\Theta$ that generalizes, we minimize an amortized loss over a distribution $\nu$ over density pairs $(\rho_1,\rho_2)$:
\begin{equation}\label{eq:amortized_loss}
\min_\Theta\ \E_{(\rho_1,\rho_2)\sim\nu}\bigl[\,\lambda_{\mathrm{mix}}\,\Ecal_{\mathrm{mix}}(\rho_1,\rho_2) + \lambda_{\mathrm{kin}}\,\Ecal_{\mathrm{kin}}(\rho_1,\rho_2)\,\bigr],
\end{equation}
with
\begin{equation}\label{eq:E_def}
\Ecal_{\mathrm{mix}}(\rho_1,\rho_2) = \|\rho(T,\cdot)\|_{\Hmone}^{2}, \quad \Ecal_{\mathrm{kin}}(\rho_1,\rho_2) = \int_0^{T}\int_{\Omega}|\bv(t,x)|^{2}\,dx\,dt,
\end{equation}
where $\rho(t,\cdot)$ solves the transport equation \eqref{eq:transport} with initial data $\rho_0=\mathbf{1}_{A}-\mathbf{1}_{A^{c}}$, $\bX(t;x)$ is the flow map \eqref{eq:ode_traj}, and both are driven by the velocity $\bv$ assembled from \eqref{eq:vel_expansion} with coefficients $\blam(t)=G_\Theta(\{x_j^{(1)}\},\{x_j^{(2)}\};t)$. The pair $(\rho_1,\rho_2)$ determines $\rho(t,\cdot)$ through the same flow map. Indeed, by linearity of \eqref{eq:transport},
\begin{equation*}
\rho(t,\cdot) \;=\; |A|\,\bX(t;\cdot)_{\#}\rho_1 \;-\; |A^{c}|\,\bX(t;\cdot)_{\#}\rho_2,
\end{equation*}
so transporting the two phases suffices. We keep the normalized measures $\rho_1,\rho_2$ as inputs because they are the laws from which the context and trajectory particles are drawn. By the volume-preserving change of variables in Section~\ref{sec:eul_to_lag}, the kinetic term in \eqref{eq:E_def} is estimated by an unweighted average over the pooled trajectory particles. Section~\ref{sec:loss_computation} gives the discrete estimators of both terms. %The weights $\lambda_{\mathrm{mix}},\lambda_{\mathrm{kin}}>0$ are fixed hyperparameters.
%=======redundant=============
%\begin{remark}\label{remark_sample2}
%For sign-changing initial data such as the two-color field $\rho_0=\mathbf 1_{A}-\mathbf 1_{A^{c}}$, the scalar $\rho_0$ cannot serve as a probability density. Training therefore draws one uniform cloud on $\Omega$ and splits it by phase: the per-phase clouds are i.i.d.\ samples of $\rho_1$ and $\rho_2$ with random sizes $n_1,n_2$ whose expectations are proportional to $|A|$ and $|A^{c}|$, and the colors enter the loss only through the signed weights of Remark~\ref{remark_sample1}. Section~\ref{sec:loss_computation} makes the resulting estimators explicit.
%\end{remark}
%=======redundant=============

\subsection{Divergence-free spectral basis on the cube}\label{sec:divfree_basis}

We construct the finite-dimensional velocity space \(V_L\) from smooth spectral modes that satisfy the incompressibility and no-penetration conditions \eqref{prop_v} exactly. Consequently, every velocity field of the form \eqref{eq:vel_expansion} satisfies these constraints for any coefficient vector \(\blam(t)\), without requiring an additional penalty.

\textbf{Two-dimensional basis.}
On $\Omega=(0,1)^{2}$, any incompressible planar field is the perpendicular gradient of a stream function $\psi$, since $\nabla\cdot(\nabla^{\perp}\psi)=0$. We take
\begin{equation*}
\psi_{mn}(x) = \sin(m\pi x_1)\sin(n\pi x_2), \quad m,n\ge 1,
\end{equation*}
and define the normalized basis
\begin{equation}\label{eq:basis_2d}
\be_{mn}(x) = C_{mn}\,\nabla^{\perp}\psi_{mn}(x) = C_{mn}\begin{pmatrix} n\pi\sin(m\pi x_1)\cos(n\pi x_2)\\[2pt] -m\pi\cos(m\pi x_1)\sin(n\pi x_2)\end{pmatrix},
\end{equation}
with $C_{mn}=2/(\pi\sqrt{m^{2}+n^{2}})$.

\textbf{Three-dimensional basis.}
On $\Omega=(0,1)^{3}$, any incompressible field is the curl of a vector potential. We use the scalar potential $\psi_{mnl}(x)=\sin(m\pi x_1)\sin(n\pi x_2)\sin(l\pi x_3)$ and take its curl around each Cartesian axis $\mathrm{e}_k$, $k\in\{1,2,3\}$. This yields three families of basis vectors, each indexed by $(m,n,l)\in\{1,\ldots,M_{\max}\}^{3}$:
\begin{equation}\label{eq:basis_3d}
\begin{aligned}
\be^{(1)}_{mnl}(x) &= C^{(1)}_{mnl}\,\nabla\times(\psi_{mnl}\,\mathrm{e}_1) = C^{(1)}_{mnl}\begin{pmatrix} 0\\ l\pi\sin(m\pi x_1)\sin(n\pi x_2)\cos(l\pi x_3)\\ -n\pi\sin(m\pi x_1)\cos(n\pi x_2)\sin(l\pi x_3)\end{pmatrix},\\
\be^{(2)}_{mnl}(x) &= C^{(2)}_{mnl}\,\nabla\times(\psi_{mnl}\,\mathrm{e}_2) = C^{(2)}_{mnl}\begin{pmatrix} -l\pi\sin(m\pi x_1)\sin(n\pi x_2)\cos(l\pi x_3)\\ 0\\ m\pi\cos(m\pi x_1)\sin(n\pi x_2)\sin(l\pi x_3)\end{pmatrix},\\
\be^{(3)}_{mnl}(x) &= C^{(3)}_{mnl}\,\nabla\times(\psi_{mnl}\,\mathrm{e}_3) = C^{(3)}_{mnl}\begin{pmatrix} n\pi\sin(m\pi x_1)\cos(n\pi x_2)\sin(l\pi x_3)\\ -m\pi\cos(m\pi x_1)\sin(n\pi x_2)\sin(l\pi x_3)\\ 0 \end{pmatrix},
\end{aligned}
\end{equation}
with the unit-norm constants $C^{(1)}_{mnl} = 2\sqrt{2}/(\pi\sqrt{n^{2}+l^{2}})$, $C^{(2)}_{mnl} = 2\sqrt{2}/(\pi\sqrt{m^{2}+l^{2}})$, and $C^{(3)}_{mnl} = 2\sqrt{2}/(\pi\sqrt{m^{2}+n^{2}})$.

It is straightforward to verify the following properties.
\begin{prop}\label{prop:basis}
The basis elements defined by \eqref{eq:basis_2d} (for $d=2$) and \eqref{eq:basis_3d} (for $d=3$) satisfy, for all admissible indices,
\begin{enumerate}[(i)]
\item $\nabla\cdot\be_\ell = 0$ in $\Omega$;
\item $\be_\ell\cdot\bn = 0$ on $\partial\Omega$;
\item $\|\be_\ell\|_{L^{2}(\Omega)} = 1$.
\end{enumerate}
\end{prop}
As a result, for  $\blam\in L^{2}(0,T;\R^{L})$ we have $\bv(t,\cdot)=\sum_{\ell=1}^{L}\lambda_\ell(t)\be_\ell \in\mathcal{F}$.
%\begin{proof}
%The first two properties follow from the identities
%\(\nabla\cdot\nabla^\perp\psi=0\) and
%\(\nabla\cdot(\nabla\times\bA)=0\), together with the vanishing sine factors on the corresponding boundary faces. The normalization follows from
%
%$$
%\int_0^1\sin^2(k\pi x)\,dx
%=
%\int_0^1\cos^2(k\pi x)\,dx
%=
%\frac12.
%$$
%
%\end{proof}
Note that the two-dimensional modes are mutually \(L^2\)-orthonormal, whereas the three-dimensional modes are individually normalized but not mutually orthogonal.
Accordingly, for the 3D modes,
$$
\int_\Omega|\bv(t,x)|^2\,dx
=
\blam(t)^\top G\blam(t),
\qquad
G_{\ell k}:=
\langle\be_\ell,\be_k\rangle_{L^2(\Omega)}.
$$
%If desired, the retained modes may be orthonormalized once with respect to the \(L^2\) inner product; this preserves their smoothness, divergence-free property, and no-penetration condition.

\begin{remark}[Truncation, basis size, and one omission]\label{rem:Mmax}
Truncating each index by $M_{\max}$ gives $L=M_{\max}^{2}$ modes in 2D and $L=3M_{\max}^{3}$ in 3D. We use $M_{\max}=2$, giving $L=4$ and $L=24$, respectively. In 3D, a zero frequency index makes the scalar potential vanish, so all retained modes have strictly positive indices. The resulting velocity space contains no nonzero flow that is exactly invariant along a coordinate direction.
\end{remark}

\subsection{Transformer architecture}\label{sec:context_encoder}

\begin{figure}[!ht]
\centering
\begin{tikzpicture}[
    node distance=0.55cm and 2.2cm,
    block/.style={rectangle, draw, rounded corners=2pt, minimum width=2.5cm, minimum height=0.55cm, align=center, fill=blue!6, font=\footnotesize},
    smallblock/.style={rectangle, draw, rounded corners=2pt, minimum width=2.2cm, minimum height=0.48cm, align=center, fill=green!6, font=\footnotesize},
    resultblock/.style={rectangle, draw, rounded corners=2pt, minimum width=2.6cm, minimum height=0.55cm, align=center, fill=orange!8, font=\footnotesize},
    arrow/.style={-{Stealth[length=4.5pt]}, thick},
    dasharrow/.style={-{Stealth[length=4.5pt]}, thick, dashed, blue!60},
    stagelab/.style={font=\scriptsize\itshape, text=gray!60!black},
]
\node[block] (rho1) {$\rho_1$ points $\{x_j^{(1)}\}$};
\node[block, right=2.4cm of rho1] (rho2) {$\rho_2$ points $\{x_j^{(2)}\}$};
\node[smallblock, below=0.5cm of rho1] (lift1) {lift $+\,g_1$};
\node[smallblock, below=0.5cm of rho2] (lift2) {lift $+\,g_2$};
\node[smallblock, below=0.4cm of lift1] (sa1) {self-attn};
\node[smallblock, below=0.4cm of lift2] (sa2) {self-attn};
\node[smallblock, below=0.55cm of sa1] (ca1) {cross-attn $\rho_1\leftarrow\rho_2$};
\node[smallblock, below=0.55cm of sa2] (ca2) {cross-attn $\rho_2\leftarrow\rho_1$};
\node[smallblock, below=0.7cm of $(ca1)!0.5!(ca2)$] (joint) {concat $+$ joint self-attn};
\node[smallblock, below=0.4cm of joint] (pool) {masked mean pool $+$ LN};
\node[resultblock, below=0.4cm of pool] (zvec) {context $z\in\R^{h}$};
\node[smallblock, below=0.55cm of zvec] (head) {$[z\,;\,e_{\mathrm{time}}(t)]\to\mathrm{MLP}\to\lambda_{\max}\tanh(\cdot)$};
\node[resultblock, below=0.4cm of head] (lam) {$\blam(t)\in\R^{L}$};
\node[resultblock, below=0.4cm of lam] (vel) {$\bv(t,x)=\sum_{\ell=1}^{L}\lambda_\ell(t)\,\be_\ell(x)$};
\draw[arrow] (rho1) -- (lift1);
\draw[arrow] (rho2) -- (lift2);
\draw[arrow] (lift1) -- (sa1);
\draw[arrow] (lift2) -- (sa2);
\draw[arrow] (sa1) -- (ca1);
\draw[arrow] (sa2) -- (ca2);
\draw[dasharrow] (sa2.south west) -- (ca1.north east);
\draw[dasharrow] (sa1.south east) -- (ca2.north west);
\draw[arrow] (ca1.south) -- ++(0,-0.25) -| (joint.north west);
\draw[arrow] (ca2.south) -- ++(0,-0.25) -| (joint.north east);
\draw[arrow] (joint) -- (pool);
\draw[arrow] (pool) -- (zvec);
\draw[arrow] (zvec) -- (head);
\draw[arrow] (head) -- (lam);
\draw[arrow] (lam) -- (vel);
\node[stagelab, left=0.50cm of sa1, rotate=90, anchor=south] {Stage 1: encoder (once)};
\node[stagelab, left=1.1cm of head, rotate=90, anchor=south] {Stage 2: head (per step)};
\end{tikzpicture}
\caption{Architecture of the in-context mixing operator $G_\Theta$. Stage 1 encodes the two point clouds into a context vector $z$. Stage 2 maps $(z,t)$ to the basis coefficients $\blam(t)$ using the projected time feature $e_{\mathrm{time}}(t)$ defined in this subsection. Dashed arrows indicate cross-attention flow of keys and values.}
\label{fig:architecture}
\end{figure}
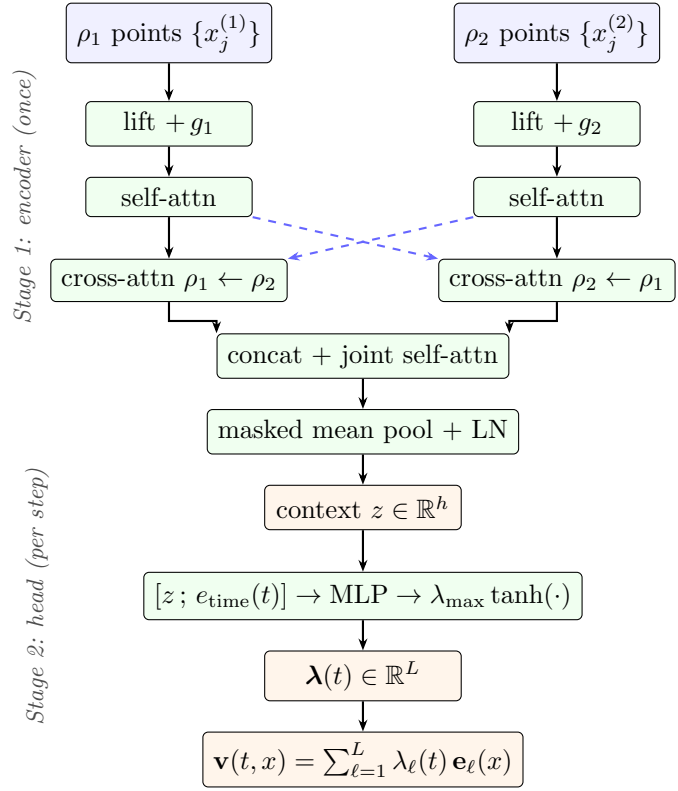

The network $G_\Theta$ is split into two stages following the cost structure of the inner loop. The \emph{context encoder} maps the two point clouds to a single vector $z\in\R^{h}$, where $h$ is the hidden dimension (Table~\ref{tab:hyperparams}), and is invoked once per trajectory. The \emph{coefficient head} maps $(z,t)$ to $\blam(t)\in\R^{L}$ and is queried at the RK4 stage times (Algorithm~\ref{alg:rk4}). The overall data flow is shown in Figure~\ref{fig:architecture}.

\textbf{Stage 1: two-density context encoder.} Each input point $x_j^{(k)}\in\R^{d}$ is lifted to $\R^{h}$ by a shared linear map $x\mapsto W_{\mathrm{lift}}x+b_{\mathrm{lift}}$, and a learned group embedding $g_k\in\R^{h}$ is added to distinguish the two phases. The two token sets are processed by $n_{\mathrm{self}}$ per-group layers of a pre-LayerNorm transformer encoder \cite{vaswani2017attention}. A bidirectional cross-attention block follows: group-$1$ tokens attend to group-$2$ tokens (queries from group $1$, keys and values from group $2$) and vice versa, so that the encoder models the interface between the two phases. The two sets are then concatenated into a joint sequence of length $n_1+n_2$, processed by $n_{\mathrm{joint}}$ further joint transformer layers, and reduced to a single vector $z\in\R^{h}$ by masked mean pooling followed by LayerNorm. The mean pool enforces permutation invariance within each group and sampling invariance as $n_1,n_2$ grow \cite{zaheer2017deepsets,huang2025unsupervised}.

\textbf{Stage 2: coefficient head.} The normalized query time $\tau=t/T$ is embedded by sinusoidal encoding with logarithmically spaced frequencies,
\begin{equation*}
\phi(\tau) = \big(\sin(\omega_1\tau),\ldots,\sin(\omega_{d_t/2}\tau),\cos(\omega_1\tau),\ldots,\cos(\omega_{d_t/2}\tau)\big)\in\R^{d_t},
\end{equation*}
where $d_t\in 2\N$ is the time-embedding dimension (Table~\ref{tab:hyperparams}). The resulting feature is mapped to $\R^{h}$ by the affine transformation
\begin{equation*}
e_{\mathrm{time}}(t)=W\phi(t/T)+b.
\end{equation*}
Here $W\in\R^{h\times d_t}$ and $b\in\R^{h}$ are trainable parameters included in $\Theta$. They are learned during training and shared across all query times. The time dependence of $e_{\mathrm{time}}(t)$ enters only through $\phi(t/T)$. The concatenation $[z\,;\,e_{\mathrm{time}}(t)]\in\R^{2h}$ is fed through a three-layer MLP ($2h\to 4h\to 2h\to L$) with GELU activations, and the output is passed through a bounded nonlinearity
\begin{equation*}
\blam(t) = \lambda_{\max}\,\tanh(\mathrm{MLP}([z;\,e_{\mathrm{time}}(t)])), \qquad \lambda_{\max}=0.5,
\end{equation*}
which constrains $\|\blam(t)\|_\infty\le\lambda_{\max}$ and prevents the velocity from destabilizing the ODE integrator. The output layer is zero-initialized so that the untrained network produces the trivial control $\bv\equiv 0$: at initialization particles remain stationary and the loss equals the initial mix-norm, which provides a well-conditioned starting point for optimization.

\textbf{Cost allocation.} The context encoder costs $\mathcal{O}((n_1+n_2)^2h)$ and is evaluated once per trajectory. Each coefficient-head evaluation costs $\mathcal{O}(h^2+hL)$. Algorithm~\ref{alg:rk4} uses three such evaluations per step, sharing the midpoint coefficients between the two intermediate stages.

%% ====================================================================
\section{Training methodology}\label{sec:training}

Training minimizes the amortized objective \eqref{eq:amortized_loss} by differentiating through an RK4 particle rollout and the terminal mix-norm estimator. Algorithm~\ref{alg:rk4} gives one integration step, and Algorithm~\ref{alg:train} combines the rollout, loss evaluation, and parameter update. The sampling distributions are specified in Section~\ref{sec:dataset}.

\begin{algorithm}[!ht]
\caption{One RK4 step, applied to both particle sets in parallel}\label{alg:rk4}
\small
\setlength{\abovedisplayskip}{4pt}
\setlength{\belowdisplayskip}{4pt}
\setlength{\abovedisplayshortskip}{2pt}
\setlength{\belowdisplayshortskip}{2pt}
\begin{algorithmic}[1]
\State \textbf{Input:} context $z\in\R^h$; positions $\bX^{(k)}\in[0,1]^{n_k\times d}$ at time $t$ for $k=1,2$; step $\Delta t=T/N_T$; accumulator $\Ecal_{\mathrm{kin}}$; coefficient head $\mathrm{Head}$ at the current parameters $\Theta$ and basis evaluation $E$.
\State Evaluate the coefficients at the three stage times:
\[
\begin{aligned}
\blam_1&\gets\mathrm{Head}(z,t/T),\qquad
\blam_2\gets\mathrm{Head}\bigl(z,(t+\tfrac12\Delta t)/T\bigr),\\
\blam_4&\gets\mathrm{Head}\bigl(z,(t+\Delta t)/T\bigr).
\end{aligned}
\]
\State For each phase $k=1,2$, compute the four slopes using the same coefficients:
\[
\begin{aligned}
\mathbf k_1^{(k)}&\gets E(\bX^{(k)})\,\blam_1,\\
\mathbf k_2^{(k)}&\gets E\bigl(\bX^{(k)}+\tfrac12\Delta t\,\mathbf k_1^{(k)}\bigr)\,\blam_2,\\
\mathbf k_3^{(k)}&\gets E\bigl(\bX^{(k)}+\tfrac12\Delta t\,\mathbf k_2^{(k)}\bigr)\,\blam_2,\\
\mathbf k_4^{(k)}&\gets E\bigl(\bX^{(k)}+\Delta t\,\mathbf k_3^{(k)}\bigr)\,\blam_4.
\end{aligned}
\]
\State Advance both phases and clamp the positions to $[0,1]^d$:
\[
\bX^{(k)}\gets\mathrm{clamp}\left(\bX^{(k)}+\tfrac{\Delta t}{6}\bigl(\mathbf k_1^{(k)}+2\mathbf k_2^{(k)}+2\mathbf k_3^{(k)}+\mathbf k_4^{(k)}\bigr),0,1\right).
\]
\State Accumulate the unweighted kinetic energy from the second-stage velocities:
\[
\Ecal_{\mathrm{kin}}\mathrel{+}=\frac{\Delta t}{n_1+n_2}\sum_{k=1,2}\sum_{i=1}^{n_k}|[\mathbf k_2^{(k)}]_i|^2.
\]
\State \Return $\bX^{(1)},\bX^{(2)},\Ecal_{\mathrm{kin}}$.
\end{algorithmic}
\end{algorithm}

\subsection{ODE rollout and loss}\label{sec:loss_computation}
For each batch element, we advance both phase clouds through $N_T$ steps of Algorithm~\ref{alg:rk4}, with $\Delta t=T/N_T$. At each stage, $E(\bX^{(k)})\in\R^{n_k\times L\times d}$ contains the basis values at the current particle positions, and contraction with $\blam$ gives the velocity. After each complete step, positions are clamped to $[0,1]^d$. By Proposition~\ref{prop:basis} the continuous-time flow is tangential to $\partial\Omega$, so particles do not leave the domain analytically. The clamp compensates for round-off drift and was observed to activate on fewer than $10^{-3}$ of particle-steps throughout training.

After the rollout, we evaluate the terminal mix-norm by a Fourier cosine expansion on the unit square or cube. The $L^{2}$-orthonormal eigenfunctions of the Neumann Laplacian are
\begin{equation*}
u_k(x) = \Big(\prod_{i=1}^{d}a_{k_i}\Big)\prod_{i=1}^{d}\cos(k_i\pi x_i),
\end{equation*}
where $x=(x_1,\ldots,x_d)$, $k=(k_1,\ldots,k_d)\in\mathbb{N}_0^{d}$, $a_0=1$, and $a_n=\sqrt{2}$ for $n\ge 1$. The corresponding eigenvalues are $\lambda_k=\pi^{2}|k|^{2}$, where $|k|^{2}=k_1^{2}+\cdots+k_d^{2}$ \cite{Grebenkov2013Geometrical}. The kernel integral \eqref{eq:Hminus1_kernel} therefore has the spectral representation
\begin{equation*}
\|\rho(T,\cdot)\|^{2}_{\Hmone} = \sum_{k\in\mathbb{N}_0^{d}\setminus\{0\}} \frac{|\langle\rho(T,\cdot),u_k\rangle_{L^{2}}|^{2}}{\pi^{2}|k|^{2}}.
\end{equation*}
We approximate the Fourier coefficients from the terminal particle positions and retain the nonzero modes with $|k|_\infty\le K_{\mathrm{mix}}$. If the clouds $\{\bX^{(1)}_{i}(T)\}_{i=1}^{n_1}$ and $\{\bX^{(2)}_{j}(T)\}_{j=1}^{n_2}$ are i.i.d.\ samples of the two transported phases, each Fourier coefficient $\langle\rho(T,\cdot),u_k\rangle_{L^{2}}$ is estimated without bias by the phase-weighted average
\begin{equation}\label{eq:dye_weighted}
\widehat{\rho}_k \;=\; \frac{|A|}{n_1}\sum_{i=1}^{n_1}u_k\bigl(\bX^{(1)}_{i}(T)\bigr)\;-\;\frac{|A^{c}|}{n_2}\sum_{j=1}^{n_2}u_k\bigl(\bX^{(2)}_{j}(T)\bigr),
\end{equation}
where the phase volumes $|A|$ and $|A^{c}|$ compensate for the normalization of $\rho_1$ and $\rho_2$. In our implementation the two clouds are not drawn separately. A single uniform sample of size $N=n_1+n_2$ is split by phase (Section~\ref{sec:dataset}), so the empirical fractions $n_1/N$ and $n_2/N$ are consistent estimators of $|A|$ and $|A^{c}|$, and substituting them into \eqref{eq:dye_weighted} yields the signed average $\frac1N\bigl(\sum_i u_k(\bX^{(1)}_{i}(T))-\sum_j u_k(\bX^{(2)}_{j}(T))\bigr)$. Our uniform sampling is therefore consistent with the weighted format \eqref{eq:dye_weighted}. Finally, the weights are centered as in \eqref{eq:lagrangian_is}: writing $w_i=\pm1$ for the color of particle $i$ in the pooled cloud $\{\bX_i(T)\}_{i=1}^{N}$ and $\tilde w_i=w_i-\frac1N\sum_{j=1}^{N}w_j$, the training estimator is
\begin{equation}\label{eq:dye_cost}
\widehat{\rho}_k = \frac1N\sum_{i=1}^{N}\tilde w_i\,u_k\bigl(\bX_i(T)\bigr),
\qquad
\Ecal_{\mathrm{mix}} = \sum_{\substack{k\in\{0,\ldots,K_{\mathrm{mix}}\}^{d}\\ k\ne 0}}\frac{|\widehat{\rho}_k|^{2}}{\pi^{2}|k|^{2}},
\end{equation}
which approximates the squared mix-norm of $P_0\rho(T,\cdot)$ by the truncated spectral sum. This is the estimator used in training and Lagrangian evaluation. Eulerian verification in Section~\ref{sec:verification} instead uses the elliptic formulation \eqref{eq:Hminus1_identity}. Centering matters only for unbalanced phases: when $n_1=n_2$ it reduces to the uncentered signed average.

With the midpoint quadrature of Algorithm~\ref{alg:rk4}, the full per-sample training loss is
\begin{equation}\label{eq:total_loss}
\Lcal(\Theta;\rho_1,\rho_2)
=\lambda_{\mathrm{mix}}\Ecal_{\mathrm{mix}}
+\lambda_{\mathrm{kin}}\frac{\Delta t}{N}
\sum_{i=1}^{N}\sum_{s=0}^{N_T-1}
\bigl|\bv\bigl(t_{s+1/2},\bX_i(t_{s+1/2})\bigr)\bigr|^2,
\end{equation}
where $\Ecal_{\mathrm{mix}}$ is given by \eqref{eq:dye_cost}, $t_{s+1/2}=(s+\tfrac12)\Delta t$, and $\bX_i(t_{s+1/2})$ is the Euler predictor used to evaluate $\mathbf k_2^{(k)}$ in Algorithm~\ref{alg:rk4}. Every quantity in \eqref{eq:total_loss} is a smooth function of the network parameters $\Theta$.
Gradients of \eqref{eq:total_loss} with respect to $\Theta$ are computed by backpropagation through the full RK4 graph. This discretize-then-optimize approach yields exact gradients of the discrete loss and avoids the implementation cost of an adjoint ODE.

\begin{remark}[Kernel singularity and spectral smoothing]\label{rem:kernel_smoothing}
The spectral sum in \eqref{eq:dye_cost} uses the truncated Neumann Green kernel
\begin{equation*}
\mathcal{K}_K(x,y):=\sum_{\substack{k\in\{0,\ldots,K_{\mathrm{mix}}\}^{d}\\ k\ne 0}}\frac{u_k(x)\,u_k(y)}{\pi^{2}|k|^{2}}.
\end{equation*}
The exact kernel $\mathcal K$ is singular on the diagonal, so its self-interaction terms are undefined for an atomic particle measure. The truncated kernel $\mathcal K_K$ is smooth and bounded, making \eqref{eq:dye_cost} well defined and differentiable in the particle positions, including the diagonal terms. The cutoff resolves length scales of order $1/K_{\mathrm{mix}}$, while the diagonal contribution is of order $\mathcal K_K(x,x)/N$. Since $\mathcal K_K(x,x)$ grows like $\log K_{\mathrm{mix}}$ in 2D and like $K_{\mathrm{mix}}$ in 3D, increasing the cutoff balances finer spectral resolution against particle noise at fixed $N$.
\end{remark}

\subsection{Optimization}\label{sec:optimization}
We use Adam with base learning rate $\alpha$, warmed up linearly from $0.01\alpha$ to $\alpha$ over the first $N_{\mathrm{iter}}/50$ iterations and then cosine-decayed to $0.01\alpha$ over the remainder. Gradients are clipped in norm at $\gamma/\sqrt{L}$. The $1/\sqrt{L}$ scaling normalizes the update magnitude per coefficient and stabilizes training across different basis sizes. Every $N_{\mathrm{eval}}$ iterations the network is evaluated on a fresh batch of density pairs drawn from both the train and test parameter distributions (Section~\ref{sec:dataset}). The mix-norm and kinetic-energy components are logged separately. The full training loop is Algorithm~\ref{alg:train}.

\begin{algorithm}[!ht]
\caption{Amortized training of the mixing operator}\label{alg:train}
\begin{algorithmic}[1]
\State \textbf{Input:} partition family with parameter domain $\mathcal T$ and test hole $H$ (Section~\ref{sec:dataset}); basis $\{\be_\ell\}_{\ell=1}^L$ (Section~\ref{sec:divfree_basis}); hyperparameters in Table~\ref{tab:hyperparams}; evaluation period $N_{\mathrm{eval}}$.
\State \textbf{Initialize:} network parameters $\Theta$ as in Section~\ref{sec:context_encoder}, with zero output-layer weights and biases; Adam and its learning-rate schedule as in Section~\ref{sec:optimization}; step size $\Delta t=T/N_T$.
\For{$i=1,\ldots,N_{\mathrm{iter}}$}
  \State \textbf{Batch draw.} Sample $B$ parameters $p_b\sim\mathrm{Uniform}(\mathcal T\setminus H)$. For each $b=1,\ldots,B$, construct independent context and trajectory clouds by uniform sampling and phase classification as in Section~\ref{sec:dataset}.
  \State \textbf{Context encoding.} Compute $z_b\gets\mathrm{Encoder}(\{x_j^{(1,b)}\},\{x_j^{(2,b)}\})$ once for each batch element.
  \State \textbf{RK4 rollout.} Set $\Ecal_{\mathrm{kin}}^{(b)}\gets0$ and apply Algorithm~\ref{alg:rk4} for $N_T$ steps to obtain the terminal clouds $\bX_T^{(k,b)}$ and kinetic energy $\Ecal_{\mathrm{kin}}^{(b)}$.
  \State \textbf{Terminal mix-norm.} Evaluate $\Ecal_{\mathrm{mix}}^{(b)}$ from the terminal clouds using \eqref{eq:dye_cost}.
  \State \textbf{Loss and gradient.} Form the batch loss
    \[
      \Lcal\gets\frac1B\sum_{b=1}^B\bigl(\lambda_{\mathrm{mix}}\Ecal_{\mathrm{mix}}^{(b)}+\lambda_{\mathrm{kin}}\Ecal_{\mathrm{kin}}^{(b)}\bigr),
    \]
    and compute $\nabla_\Theta\Lcal$ by backpropagation through the RK4 graph.
  \State \textbf{Parameter update.} Apply Adam to the gradient clipped at $\gamma/\sqrt L$ and advance the learning-rate schedule.
  \State \textbf{Evaluation.} Every $N_{\mathrm{eval}}$ iterations, evaluate and log the mix-norm and kinetic-energy components on fresh batches from $\mathcal T\setminus H$ and $H$.
\EndFor
\State \Return $\Theta$.
\end{algorithmic}
\end{algorithm}

%% ====================================================================
\subsection{Parametric partition families}\label{sec:dataset}

The amortized objective of Section~\ref{sec:solution_op} takes expectation over a distribution $\nu$ of density pairs $(\rho_1,\rho_2)$. This section specifies $\nu$ through eight parametric families of two-phase partitions of the square or cube, four in 2D and four in 3D, with intrinsic dimensions $\dim\in\{1,2,3,4,6,9\}$. The parameter domain $\mathcal{T}$ and test hole $H$ of every family are listed in Table~\ref{tab:families}. The following paragraphs introduce only the predicate $\chi_p$.

\textbf{Common construction.} Each family is indexed by a parameter $p\in\mathcal{T}\subset\R^{\dim}$ and a predicate $\chi_p:\Omega\to\{0,1\}$ encoding the partition. Given $p$, the associated density pair consists of the uniform measures on the two complementary sets
\begin{equation}\label{eq:partition_pair}
\rho_1^{(p)} = \mathrm{Uniform}\bigl(\{x\in\Omega:\chi_p(x)=1\}\bigr),\qquad \rho_2^{(p)} = \mathrm{Uniform}\bigl(\{x\in\Omega:\chi_p(x)=0\}\bigr).
\end{equation}
The parameter domain $\mathcal{T}$ is equipped with a test hole $H\subset\mathcal{T}$. Training parameters are drawn uniformly on $\mathcal{T}\setminus H$ by rejection sampling, and test parameters are drawn uniformly on $H$. This hole-out split measures interpolation generalization across a contiguous region of the parameter manifold rather than through a random shuffle: the test parameters form a contiguous block $H$ that the training set excludes. Throughout the section we write $x_{\ast}=(\tfrac12,\ldots,\tfrac12)^{\top}\in\Omega$ for the center of the square or cube: the domain is $\Omega=[0,1]^{2}$ for the four 2D families below and $\Omega=[0,1]^{3}$ for the four 3D families. For the two Fourier level-set families (2D-4 and 3D-4), draws whose two phases are strongly imbalanced (either class occupying less than $20\%$ of $\Omega$) are rejected, so that both $\rho_1^{(p)}$ and $\rho_2^{(p)}$ are well-resolved by finite particle counts. For the other families the parameter domains $\mathcal{T}$ in Table~\ref{tab:families} keep both phases nonempty.

\textbf{Family 2D-1: rotation line ($\dim=1$).}
Let $n(\theta)=(\cos\theta,\sin\theta)^{\top}$. Family 2D-1 bisects the square by a line through the center with rotating normal,
\begin{equation*}
\chi_{\theta}(x) \;=\; \mathbf{1}\bigl[(x-x_{\ast})\cdot n(\theta)<0\bigr].
\end{equation*}
Each member is a diameter of the square rotated by $\theta$, and the two phases have equal area.

\textbf{Family 2D-2: general line ($\dim=2$).}
Adding a signed offset $b$ to the predicate of Family 2D-1 gives the two-parameter family
\begin{equation*}
\chi_{\theta,b}(x) \;=\; \mathbf{1}\bigl[(x-x_{\ast})\cdot n(\theta)<b\bigr].
\end{equation*}
Unlike Family 2D-1, the two phases have unequal area in general.

\textbf{Family 2D-3: circle ($\dim=3$).}
A disc of variable center $c=(c_1,c_2)$ and radius $r$,
\begin{equation*}
\chi_{c,r}(x) \;=\; \mathbf{1}\bigl[\|x-c\|<r\bigr].
\end{equation*}
The interface is curved and the two phases are generally unbalanced.

\textbf{Family 2D-4: Fourier level set ($\dim=2K$ with $K=3$).}
A smooth non-convex partition defined by the zero set of a truncated bivariate Fourier series with coefficients $p=(a_1,\ldots,a_K,b_1,\ldots,b_K)$,
\begin{equation}\label{eq:fourier_2d}
\begin{aligned}
	&f_{p}(x) \;=\; \sum_{k=1}^{K}\Bigl[a_k\sin(k\pi x_1)\cos(k\pi x_2)+b_k\cos(k\pi x_1)\sin(k\pi x_2)\Bigr],\\ 
	&\chi_{p}(x) \;=\; \mathbf{1}\bigl[f_{p}(x)>0\bigr].
\end{aligned}
\end{equation}
The resulting interfaces are generically curved, multiply connected, and not isotopic to a line or a circle. This family tests generalization beyond the simple affine and radial geometries of Families 2D-1 to 2D-3.

\medskip
The 3D families mirror the 2D families, with planes replacing lines, balls replacing discs, and a trivariate Fourier level set replacing the bivariate one. Let $n(\theta,\varphi)=(\sin\varphi\cos\theta,\sin\varphi\sin\theta,\cos\varphi)^{\top}$ denote the unit normal in spherical coordinates, where restricting $(\theta,\varphi)$ to $[0,\pi)\times[0,\pi)$ covers the unoriented directions (opposite normals give the same plane).

\textbf{Family 3D-1: rotation plane ($\dim=2$).}
A bisection of the cube by a plane through the center, parametrized by the two spherical angles of its normal:
\begin{equation*}
\chi_{\theta,\varphi}(x) \;=\; \mathbf{1}\bigl[(x-x_{\ast})\cdot n(\theta,\varphi)<0\bigr].
\end{equation*}
Each member is a plane through $x_{\ast}$, and the two phases have equal volume. This is the 3D analogue of Family 2D-1.

\textbf{Family 3D-2: general plane ($\dim=3$).}
Adding a signed offset $b$ to the predicate of Family 3D-1 gives
\begin{equation*}
\chi_{\theta,\varphi,b}(x) \;=\; \mathbf{1}\bigl[(x-x_{\ast})\cdot n(\theta,\varphi)<b\bigr].
\end{equation*}
This is the 3D analogue of Family 2D-2.

\textbf{Family 3D-3: sphere ($\dim=4$).}
A ball of variable center $c=(c_1,c_2,c_3)$ and radius $r$,
\begin{equation*}
\chi_{c,r}(x) \;=\; \mathbf{1}\bigl[\|x-c\|<r\bigr].
\end{equation*}
This is the 3D analogue of Family 2D-3.

\textbf{Family 3D-4: Fourier level set ($\dim=3K$ with $K=3$).}
A smooth non-convex partition defined by the sign of a truncated trivariate Fourier series. We take the natural three-term cyclic generalization of \eqref{eq:fourier_2d}: for each frequency $k=1,\ldots,K$ the term attached to coordinate $j\in\{1,2,3\}$ carries a sine in $x_j$ and cosines in the other two coordinates, giving three scalar coefficients $(a_k,b_k,c_k)$ per frequency, so that $p=(a_1,\ldots,a_K,b_1,\ldots,b_K,c_1,\ldots,c_K)$. Explicitly,
\begin{equation}\label{eq:fourier_3d}
\begin{aligned}
f_{p}(x) \;=\;& \sum_{k=1}^{K}\Bigl[\;a_k\sin(k\pi x_1)\cos(k\pi x_2)\cos(k\pi x_3) +b_k\cos(k\pi x_1)\sin(k\pi x_2)\cos(k\pi x_3) \\ 
	&+c_k\cos(k\pi x_1)\cos(k\pi x_2)\sin(k\pi x_3)\Bigr],\\
\chi_{p}(x) \;=\;& \mathbf{1}\bigl[f_{p}(x)>0\bigr].
\end{aligned}
\end{equation}
As for Family 2D-4, this family produces genuinely three-dimensional, curved, non-convex interfaces that depend on all three coordinates, and it is the most geometrically demanding of the eight families.

\begin{table}[!ht]
\centering\small
\caption{The eight parametric partition families used in the experiments, together with their parameter domains $\mathcal{T}$ and test holes $H$. Intrinsic dimension is $\dim$. Training parameters are drawn from $\mathcal{T}\setminus H$ and test parameters from $H$. The corresponding density pairs $(\rho_1^{(p)},\rho_2^{(p)})$ form the training and test datasets. The coefficient head outputs $L=4$ basis modes in 2D and $L=24$ in 3D.}\label{tab:families}
\setlength{\tabcolsep}{3pt}
\begin{tabular}{@{}lcll@{}}
\toprule
\textbf{Family} & $\dim$ & $\mathcal{T}$ & $H$ \\
\midrule
Rotation line (2D-1) & $1$ & $\theta\in[0,\pi)$ & $\theta\in[\pi/3,\,2\pi/3]$\\
General line  (2D-2) & $2$ & $\theta\in[0,\pi),\ b\in[-0.3,\,0.3]$ & $[\pi/3,\,2\pi/3]\times[-0.1,\,0.1]$\\
Circle        (2D-3) & $3$ & $c_1,c_2\in[0.25,\,0.75],\ r\in[0.15,\,0.35]$ & $[0.4,\,0.6]^{2}\times[0.2,\,0.3]$\\
Fourier level set (2D-4) & $6$ & $a_k,b_k\in[-1,1],\ K=3$ & $a_1\in[-0.3,\,0.3]$\\
\midrule
Rotation plane (3D-1) & $2$ & $\theta,\varphi\in[0,\pi)$ & $[\pi/3,\,2\pi/3]^{2}$\\
General plane  (3D-2) & $3$ & $\theta,\varphi\in[0,\pi),\ b\in[-0.3,\,0.3]$ & $[\pi/3,\,2\pi/3]^{2}\times[-0.1,\,0.1]$\\
Sphere         (3D-3) & $4$ & $c_1,c_2,c_3\in[0.40,\,0.60],\ r\in[0.42,\,0.54]$ & $[0.45,\,0.55]^{3}\times[0.48,\,0.52]$\\
Fourier level set (3D-4) & $9$ & $a_k,b_k,c_k\in[-1,1],\ K=3$ & $a_1\in[-0.3,\,0.3]$\\
\bottomrule
\end{tabular}
\end{table}

\textbf{Particle sampling at training time.}
At each iteration the sampler returns $B$ density pairs. For each pair, $n_{\mathrm{ctx}}+n_{\mathrm{traj}}$ points are drawn uniformly in $\Omega$ and classified by $\chi_p$. The first $n_{\mathrm{ctx}}$ points form the context pool and the remaining $n_{\mathrm{traj}}$ points form the trajectory pool. Each pool is then split by phase, so the per-phase counts are random and sum to the size of the pool. The trajectory pool provides the $N=n_1+n_2=n_{\mathrm{traj}}$ particles of the loss \eqref{eq:total_loss}, and the context pool provides the encoder input. The random split ensures that both the context and the trajectory are unbiased i.i.d.\ samples of the respective densities. Different draws produce different class counts, so each batch is padded to a common size and binary masks are propagated through attention and loss computations.

%% ====================================================================
\section{Numerical experiments}\label{sec:experiments}

The experiments assess generalization within each of the eight partition families and zero-shot transfer between families of the same dimension. After specifying the common setup, we present the 2D and 3D results and verify the learned controls with an independent Eulerian solver (Section~\ref{sec:verification}).

\subsection{Experimental setup}\label{sec:experiments_setup}
All training runs use the hyperparameters in Table~\ref{tab:hyperparams}. Each family is trained independently. After training, the operator is evaluated on fresh density pairs: training pairs are drawn from parameters in $\mathcal{T}\setminus H$ and test pairs from parameters in $H$. Evaluation particles are placed on a deterministic uniform lattice ($50\times 50$ in 2D and $30\times 30\times 30$ in 3D) instead of the i.i.d.\ samples used at training time. This makes the reported mix-norm curves reproducible per parameter. The squared mix-norm $\|\rho(t,\cdot)\|^{2}_{\Hmone}$ is recorded at the RK4 snapshot times. Exponential decay rates are fitted by linear regression of $\log\|\rho(t,\cdot)\|_{\Hmone}$ against $t$ over the last $75\%$ of the trajectory.

Training uses the discrete loss \eqref{eq:total_loss}. Evaluation uses the centered estimator \eqref{eq:dye_cost} on the deterministic lattice above, so all reported mix-norms follow the mean-zero convention \eqref{eq:P0_convention}.

The spectral cutoff $K_{\mathrm{mix}}$ of the mix-norm estimator was selected by the following ablation. We retrained one family per dimension with the cutoff doubled, Family 2D-3 at $K_{\mathrm{mix}}=24$ and Family 3D-3 at $K_{\mathrm{mix}}=16$, keeping all other hyperparameters fixed. Since two operators trained at different cutoffs cannot be compared through their own training losses, both were evaluated on the same held-out partitions under one common reference estimator with a cutoff larger than either training value. The reference estimator excludes the diagonal $i=j$, removing the diagonal contribution of order $\mathcal{K}_K(x,x)/N$ discussed in Remark~\ref{rem:kernel_smoothing}. The results are collected in Table~\ref{tab:kmix_ablation}. In 2D the finer cutoff reduces the terminal squared mix-norm by about $38\%$. In 3D the two cutoffs agree to within $3\%$. We keep $K_{\mathrm{mix}}=12$ (2D) and $8$ (3D) for the results reported below, so that all families within a dimension share one cutoff. The 2D mix-norms reported below therefore sit above what the finer cutoff attains. The comparison uses one run per configuration, so the $3\%$ difference in 3D is not resolved by these data.

\begin{table}[!ht]
\centering\small
\caption{Spectral cutoff ablation. Terminal squared mix-norm on held-out test partitions, averaged over $8$ parameters (2D) and $6$ parameters (3D). Both operators of a block are evaluated under one common diagonal-free reference estimator, with reference cutoff $K_{\mathrm{mix}}=48$ in 2D and $K_{\mathrm{mix}}=32$ in 3D. The ratio normalizes each block by its first row.}\label{tab:kmix_ablation}
\begin{tabular}{llcc}
\toprule
\textbf{Family} & \textbf{Trained at} & $\Ecal_{\mathrm{mix}}$ (common ref.) & \textbf{Ratio}\\
\midrule
\multirow{2}{*}{Circle (2D-3)} & $K_{\mathrm{mix}}=12$ & $5.2\cdot10^{-4}$ & $1.00$\\
                              & $K_{\mathrm{mix}}=24$ & $3.2\cdot10^{-4}$ & $0.62$\\
\midrule
\multirow{2}{*}{Sphere (3D-3)} & $K_{\mathrm{mix}}=8$  & $1.70\cdot10^{-3}$ & $1.00$\\
                               & $K_{\mathrm{mix}}=16$ & $1.75\cdot10^{-3}$ & $1.03$\\
\bottomrule
\end{tabular}
\end{table}
%========redundant===========
%The basis cutoff $M_{\max}=2$ is a compromise between expressiveness and output dimension. Increasing $M_{\max}$ scales the coefficient head's output as $M_{\max}^{2}$ in 2D and $3M_{\max}^{3}$ in 3D, and Proposition~\ref{prop:density} shows that the $M_{\max}=2$ subspace already contains the lowest-frequency cellular modes that dominate mixing on $[0,1]^{d}$. 
%========redundant===========
All experiments were performed on the HPC provided by the Erlangen National High Performance Computing Center (NHR\@ FAU) of FAU Erlangen-N\"{u}rnberg with NVIDIA A100 GPU. The code is implemented in Python using PyTorch and is publicly available at \url{https://github.com/DCN-FAU-AvH/Mixing-OP}.

\begin{table}[!ht]
\centering\small
\caption{Training hyperparameters.}\label{tab:hyperparams}
\begin{tabular}{llc}
\toprule
\textbf{Parameter} & \textbf{Symbol} & \textbf{Value}\\
\midrule
Basis cutoff & $M_{\max}$ & $2$\\
Spectral truncation & $K_{\mathrm{mix}}$ & $12$ (2D) / $8$ (3D)\\
Hidden dimension & $h$ & $256$\\
Attention heads & $n_{\mathrm{heads}}$ & $4$\\
Self-attn layers (per group / joint) & $n_{\mathrm{self}}/n_{\mathrm{joint}}$ & $1/1$ (2D) / $2/2$ (3D)\\
Time-embedding dim & $d_t$ & $64$\\
Time horizon & $T$ & $5.0$\\
RK4 steps & $N_T$ & $20$ (2D) / $30$ (3D)\\
Context/trajectory particles & $n_{\mathrm{ctx}}/n_{\mathrm{traj}}$ & $256/256$ (2D) / $1024/1024$ (3D)\\
Batch size & $B$ & $8$ (2D) / $4$ (3D)\\
Loss weights & $\lambda_{\mathrm{mix}},\lambda_{\mathrm{kin}}$ & $1.0,\,2\cdot10^{-3}$ (2D) / $1.0,\,2\cdot10^{-2}$ (3D)\\
Learning rate & $\alpha$ & $10^{-4}$\\
Gradient clip & $\gamma$ & $1.0$\\
Training iterations & $N_{\mathrm{iter}}$ & $10^{5}$ (2D) / $8\cdot 10^{4}$ (3D)\\
\bottomrule
\end{tabular}
\end{table}

\subsection{Two-dimensional case}\label{sec:experiments_2d}
\begin{figure}[!ht]
\centering
\includegraphics[width=\textwidth]{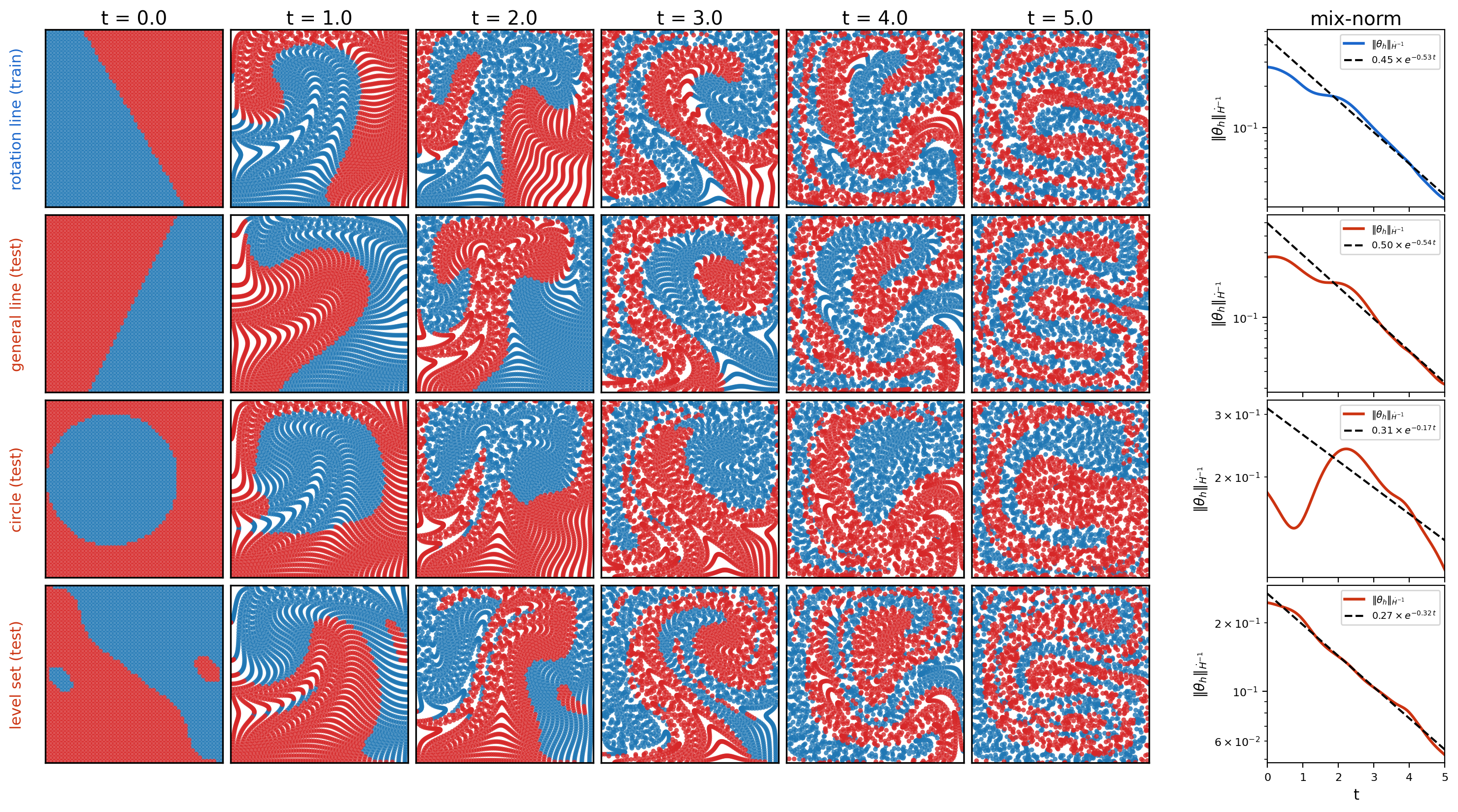}
\caption{Operator trained on 2D Rotation line ($\dim=1$).}
\label{fig:2d_rotation_line}
\end{figure}

\begin{figure}[!ht]
\centering
\includegraphics[width=\textwidth]{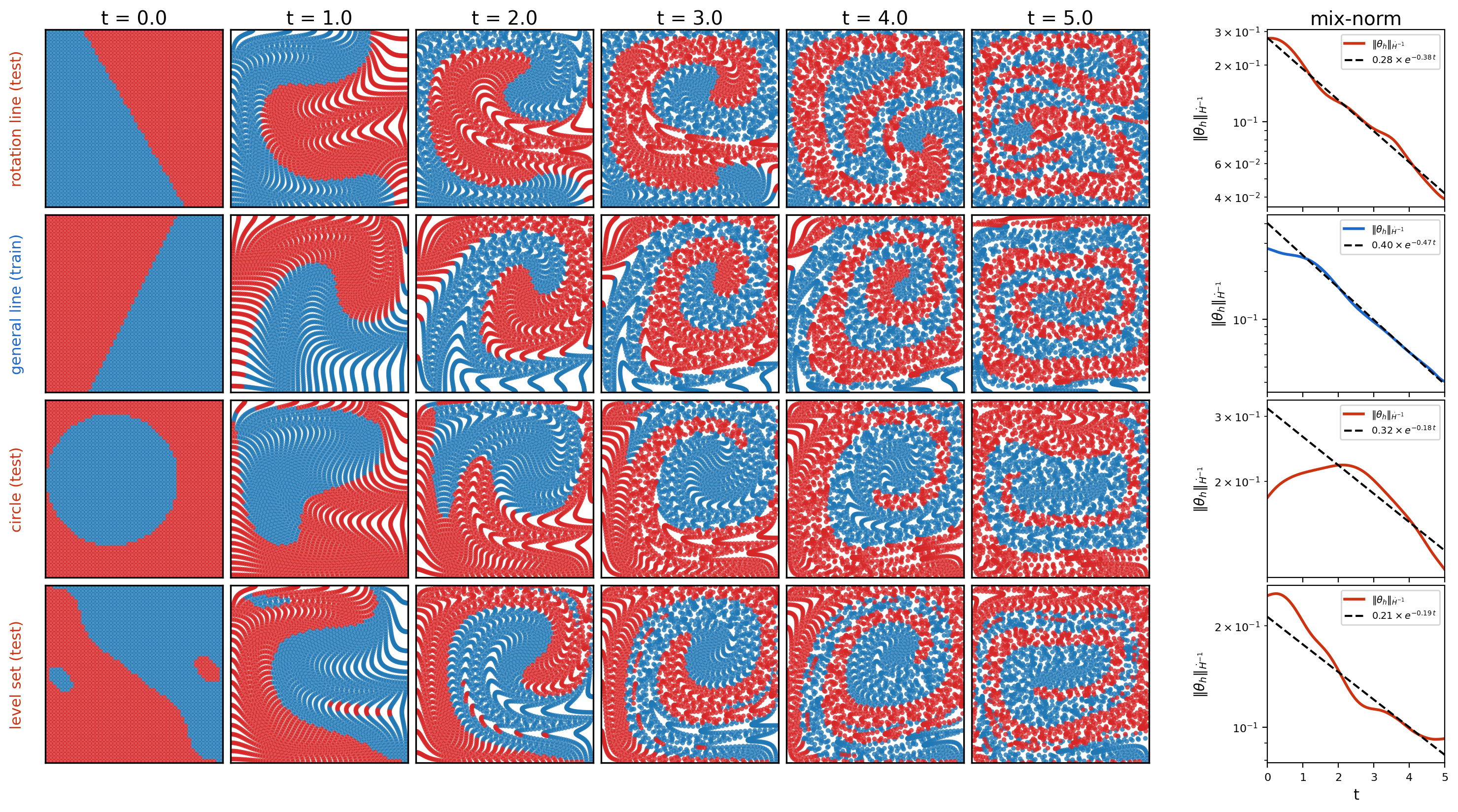}
\caption{Operator trained on 2D General line ($\dim=2$).}
\label{fig:2d_general_line}
\end{figure}

\begin{figure}[!ht]
\centering
\includegraphics[width=\textwidth]{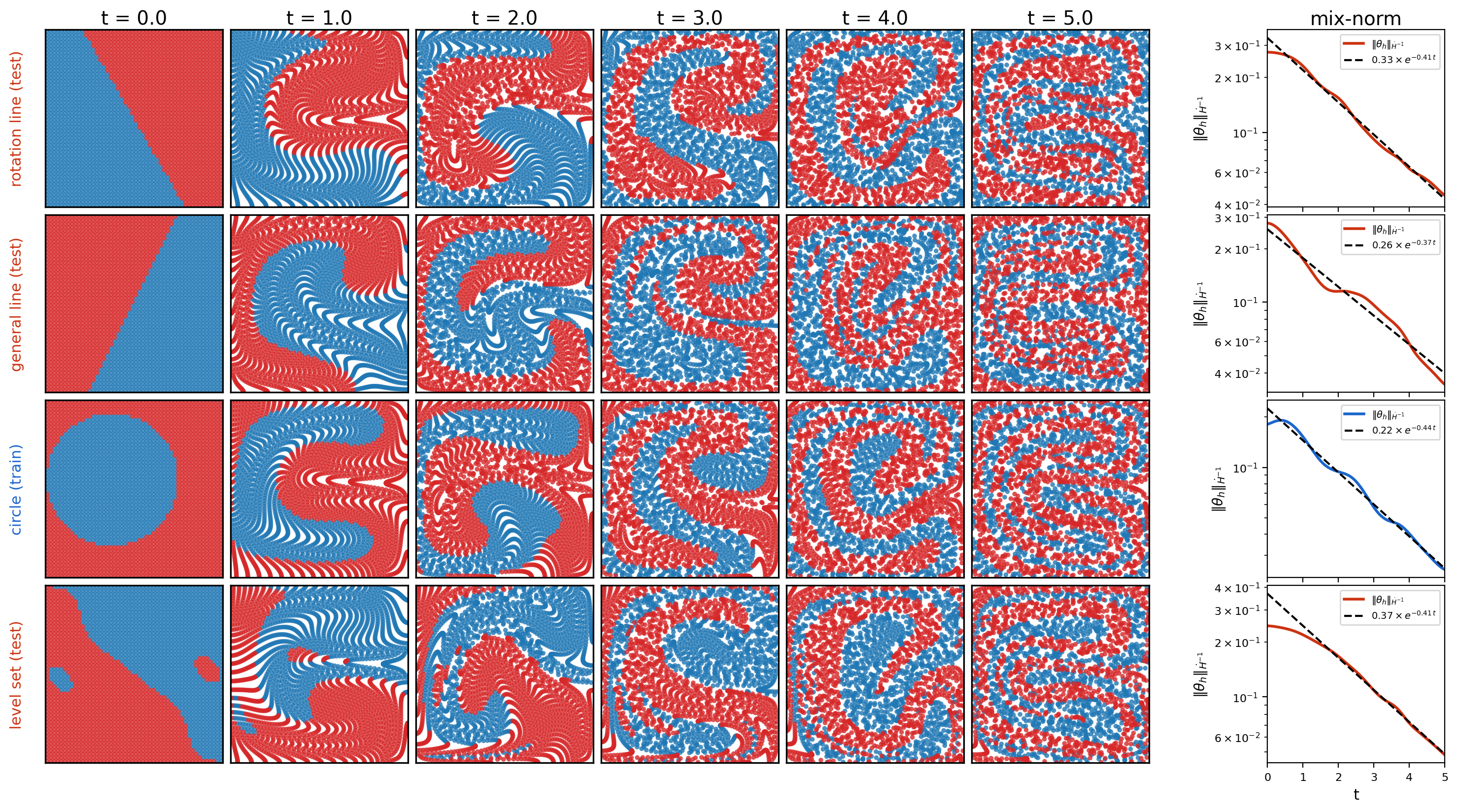}
\caption{Operator trained on 2D Circle ($\dim=3$).}
\label{fig:2d_circle}
\end{figure}

\begin{figure}[!ht]
\centering
\includegraphics[width=\textwidth]{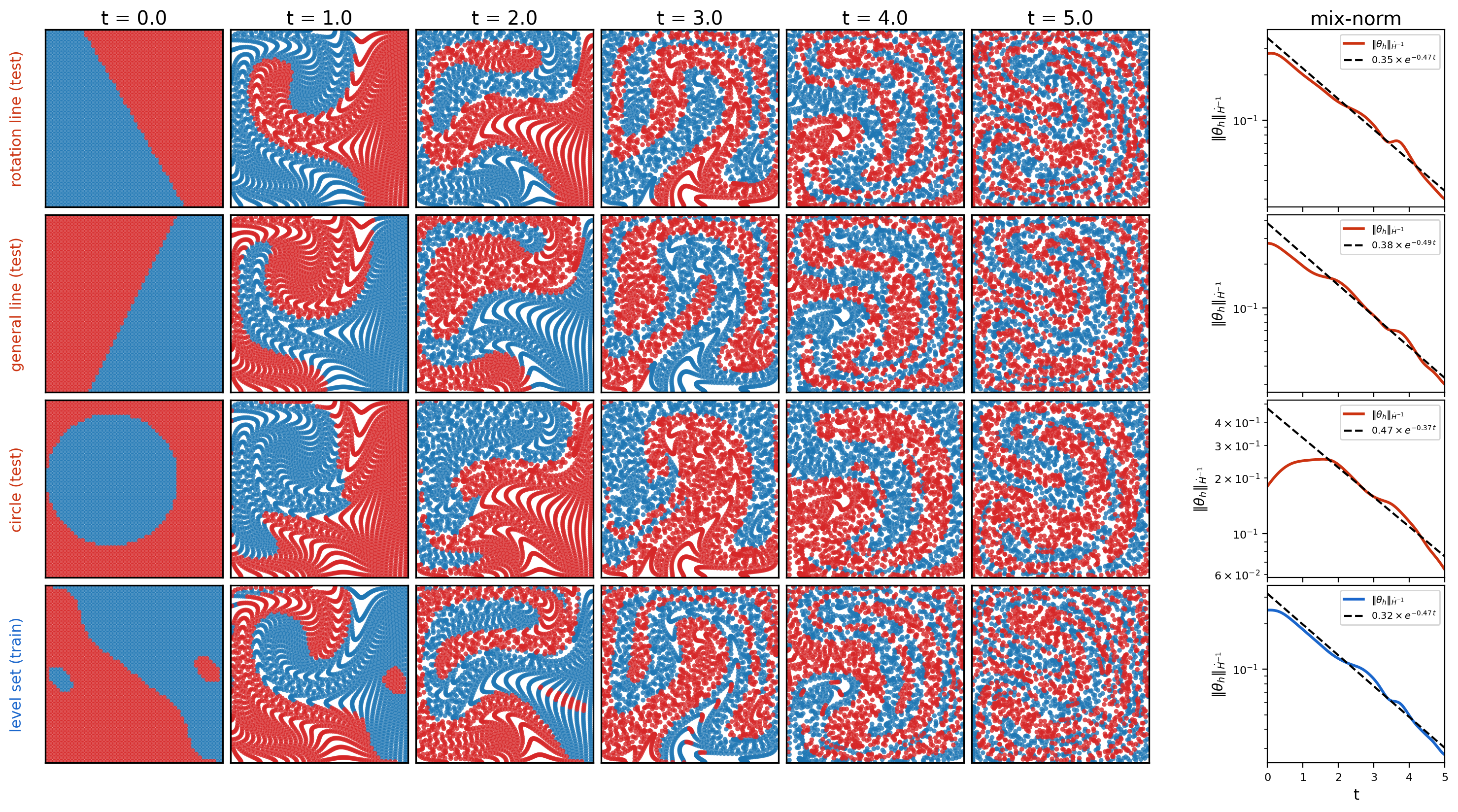}
\caption{Operator trained on 2D Fourier level set ($\dim=6$).}
\label{fig:2d_level_set}
\end{figure}

For each of the four 2D families we train an operator on that family alone and then evaluate it on all four families. Figures~\ref{fig:2d_rotation_line}--\ref{fig:2d_level_set} each contain four rows. The first row (blue spine) is the operator's native family. The remaining three rows (red spines) are the other three 2D families and measure zero-shot cross-family transfer. Every row shows six particle snapshots at $t\in\{0,T/5,\ldots,T\}$ and the $\Hmone$ mix-norm decay on a log scale with its fitted exponential $k\exp(b\,t)$.

On the native row (blue spine), the mix-norm decays approximately exponentially in all four families, and the fitted slopes $b$ have the same order of magnitude across parameter draws. The decay rate is therefore essentially parameter-independent on the family the operator is trained on. On the cross-family rows (red spines), transfer is partial: the mix-norm typically decays by about one order of magnitude over $[0,T]$, with a smaller fitted slope, and can stall when the source and target partitions have very different geometry. The learned control is also more complex as the intrinsic dimension of the native family grows: Family~2D-1 reduces to a single cellular rotation, whereas Family~2D-4 induces a filamentary pattern driven by the six-dimensional Fourier geometry.

\FloatBarrier
\subsection{Three-dimensional case}\label{sec:experiments_3d}
\begin{figure}[!ht]
\centering
\includegraphics[width=\textwidth]{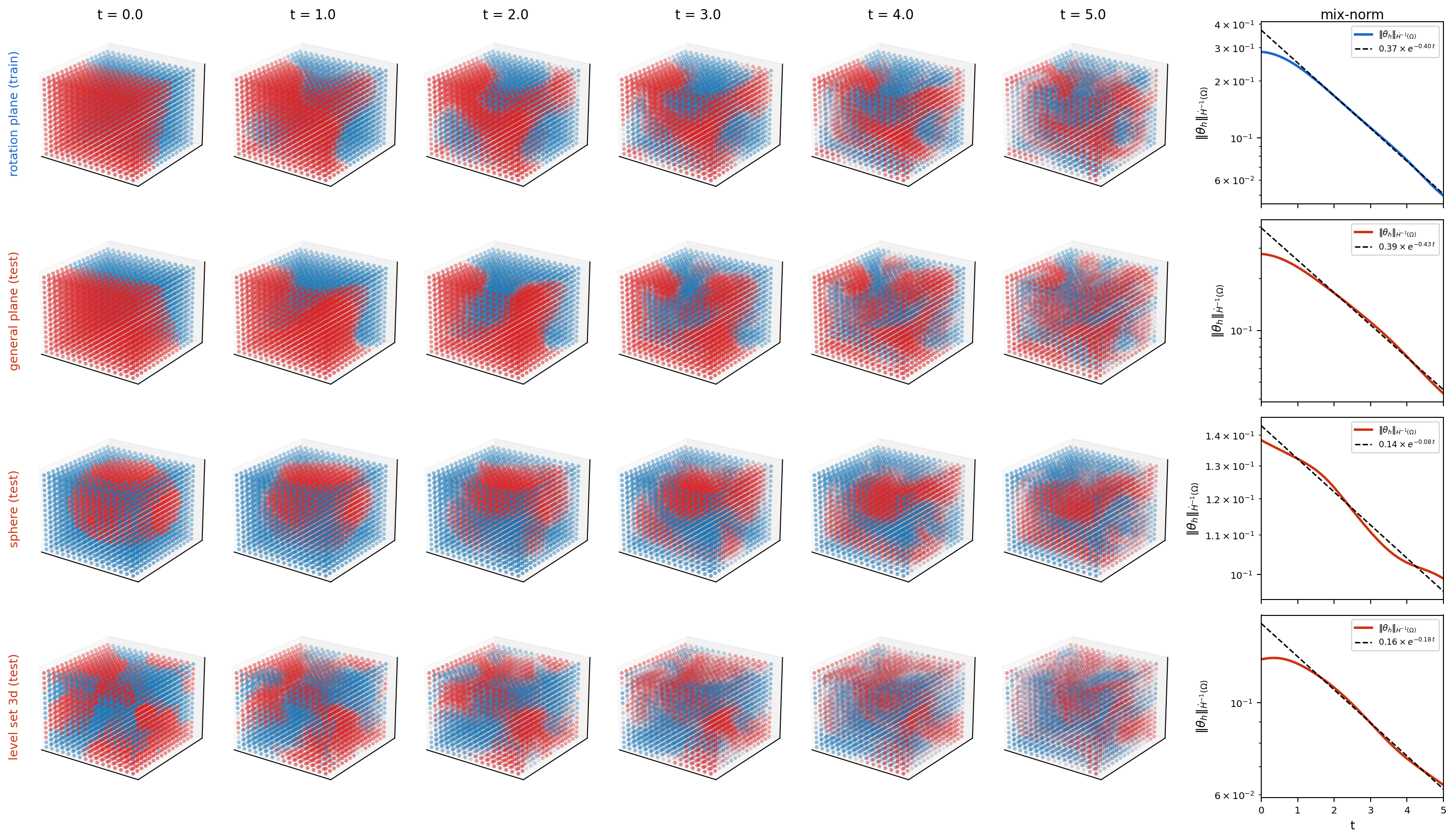}
\caption{Operator trained on 3D Rotation plane ($\dim=2$).}
\label{fig:3d_rotation_plane}
\end{figure}

\begin{figure}[!ht]
\centering
\includegraphics[width=\textwidth]{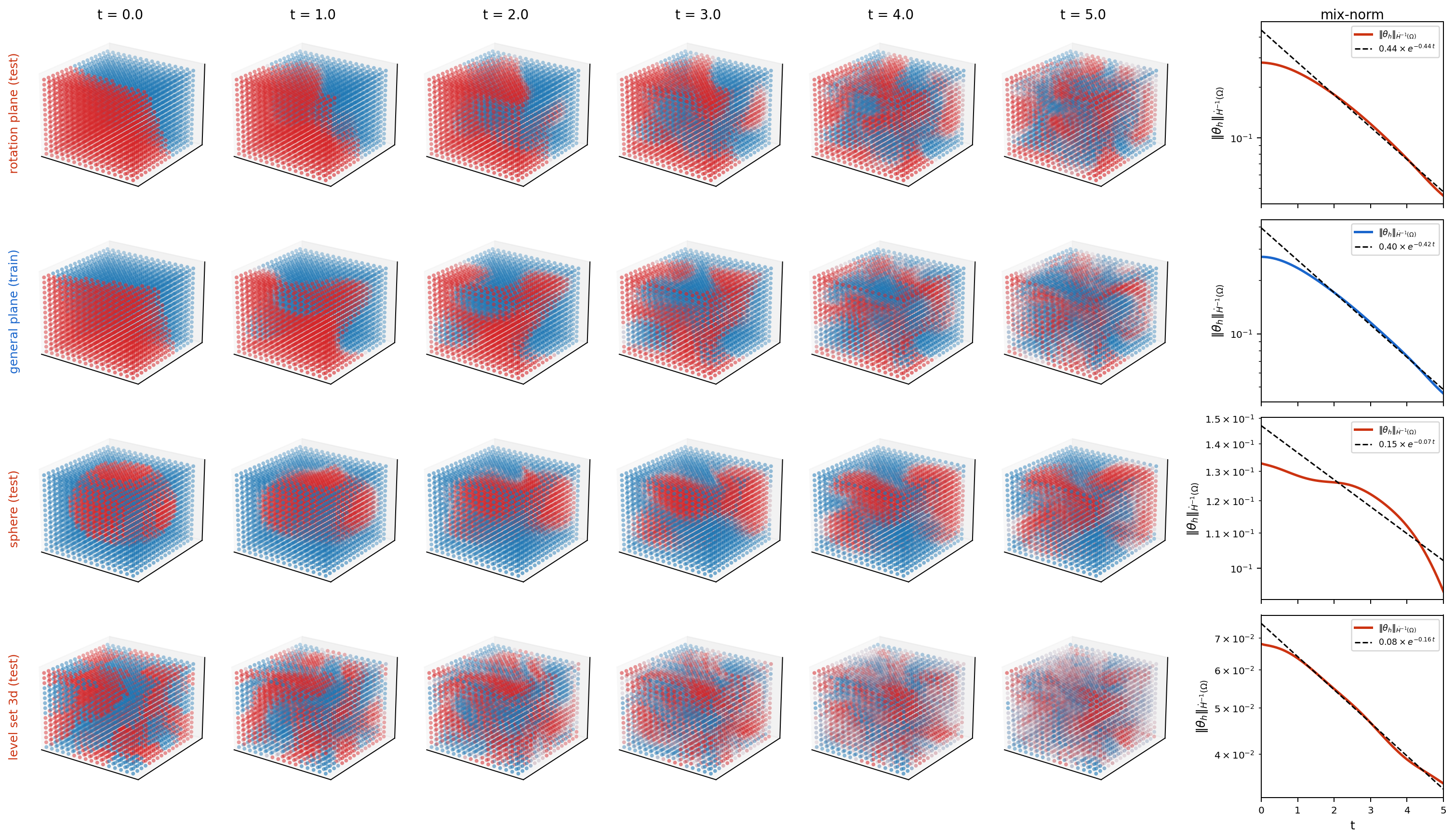}
\caption{Operator trained on 3D General plane ($\dim=3$).}
\label{fig:3d_general_plane}
\end{figure}

\begin{figure}[!ht]
\centering
\includegraphics[width=\textwidth]{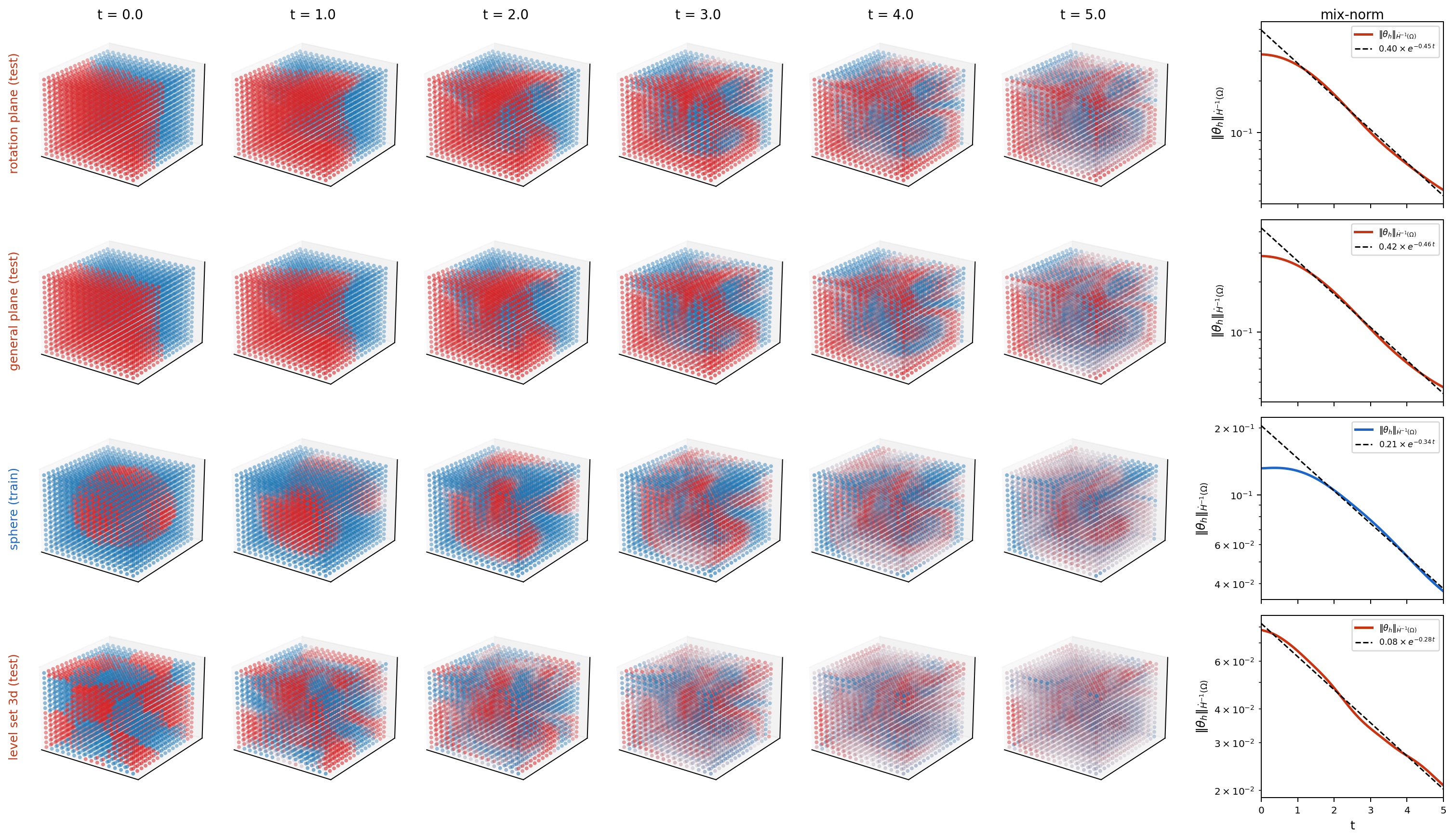}
\caption{Operator trained on 3D Sphere ($\dim=4$). }
\label{fig:3d_sphere}
\end{figure}

\begin{figure}[!ht]
\centering
\includegraphics[width=\textwidth]{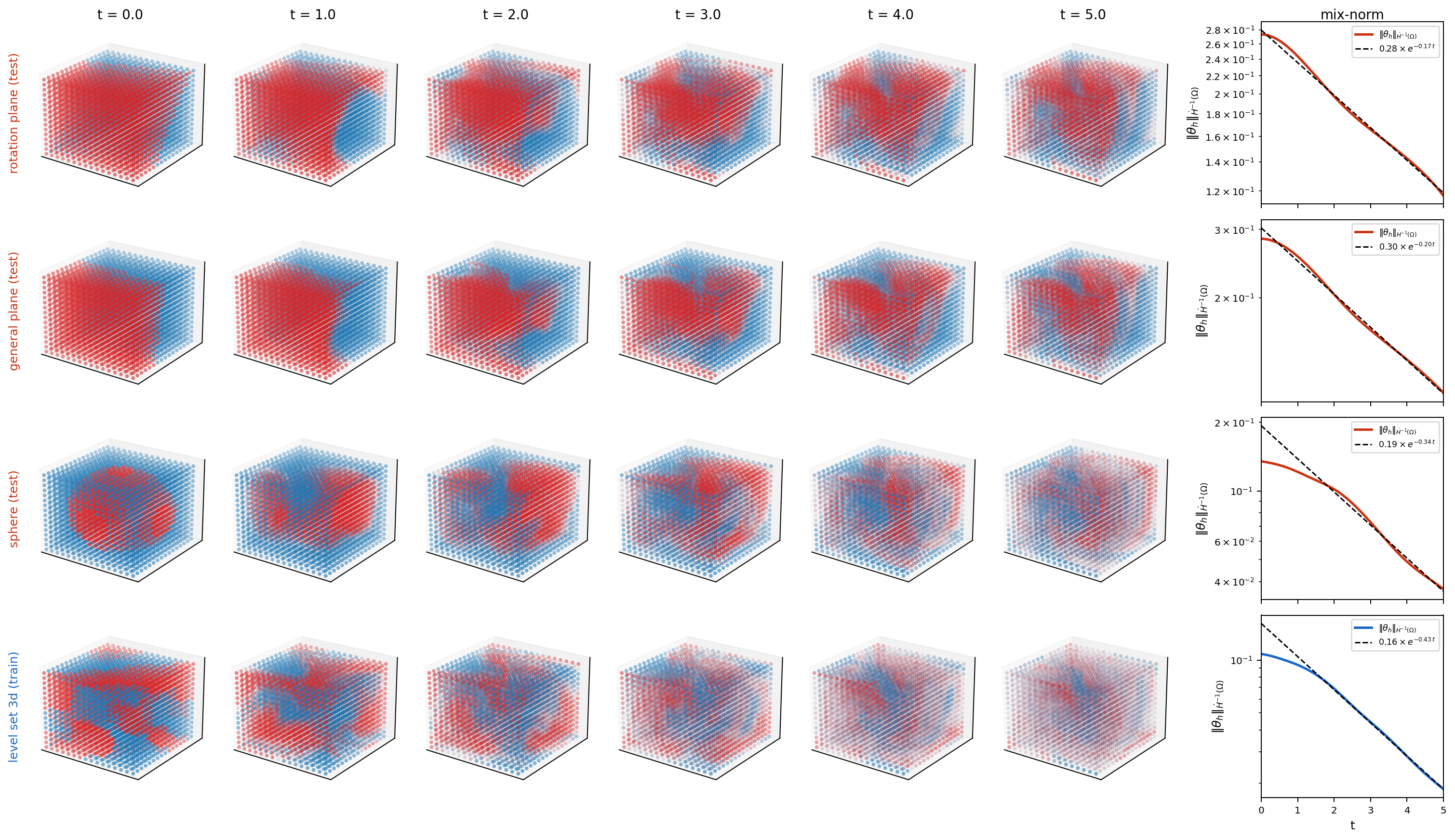}
\caption{Operator trained on 3D Level set ($\dim=9$). }
\label{fig:3d_level_set}
\end{figure}

The 3D runs use the settings in Table~\ref{tab:hyperparams}, with $L=24$ velocity modes and a $30\times30\times30$ evaluation grid. Figures~\ref{fig:3d_rotation_plane}--\ref{fig:3d_level_set} follow the same native-family and cross-family layout as in 2D.

On the native family, the 3D mix-norm again decays approximately exponentially, with fitted slopes milder than in 2D. This is consistent with the factor-of-six larger basis ($L=24$ against $L=4$ in 2D). Cross-family transfer is again partial: when the source and target partitions share broad geometric features (for instance, plane-like partitions evaluated by an operator trained on spheres) the mix-norm still decays, while on geometrically distant pairs the decay is slower or absent. Multi-family or curriculum training is a natural extension.

\FloatBarrier
\subsection{Verification against the Eulerian PDE}\label{sec:verification}
The operator is trained on a Lagrangian particle loss, but its output is a velocity field on $\Omega$. This field should produce the same mixing behaviour under the Eulerian transport equation \eqref{eq:transport}. To test this, we fix a test-parameter instance of each 2D family and read the network coefficients $\blam^{\mathrm{NN}}(t)$ at $N_T$ time midpoints. We then advance the same initial density under the same velocity field by two different discretizations: (i) an RK4 particle rollout with $6000$ particles on $[0,1]^{2}$, and (ii) a finite-volume solve of \eqref{eq:transport} on a cell-centered $N_x=256$ grid.

The finite-volume scheme follows \cite{LeVeque2002FVM}. Face velocities are reconstructed from $\blam^{\mathrm{NN}}(t)$ through the analytical integrals of the basis functions over each cell face. The scalar transport equation is advanced by a second-order MUSCL discretization with a minmod slope limiter and a Heun (RK2) time integrator, under adaptive sub-stepping that respects the CFL condition. At each snapshot time, the squared mix-norm of the gridded scalar field is obtained by solving the Neumann problem \eqref{eq:elliptic} and evaluating $\|\nabla\eta\|_{L^{2}}^{2}$ via \eqref{eq:Hminus1_identity}.

\begin{figure}[!ht]
\centering
\includegraphics[width=\textwidth]{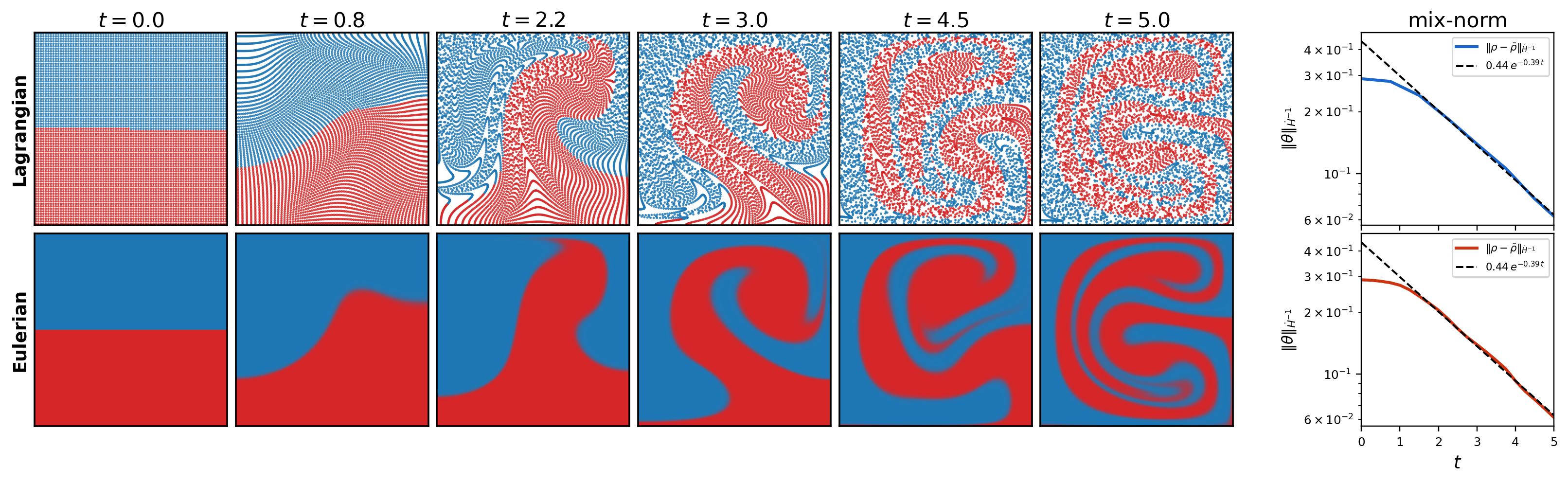}
\caption{Verification for Family 2D-1 (rotation line).}
\label{fig:eulerian_rotation_line}
\end{figure}

\begin{figure}[!ht]
\centering
\includegraphics[width=\textwidth]{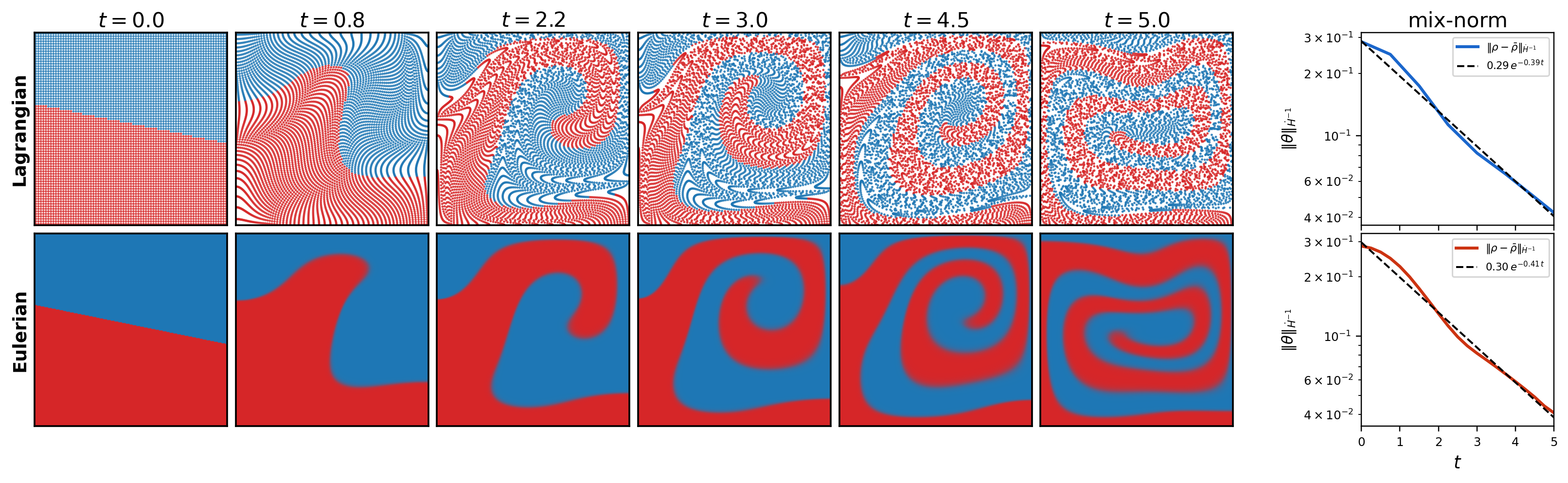}
\caption{Verification for Family 2D-2 (general line).}
\label{fig:eulerian_general_line}
\end{figure}

\begin{figure}[!ht]
\centering
\includegraphics[width=\textwidth]{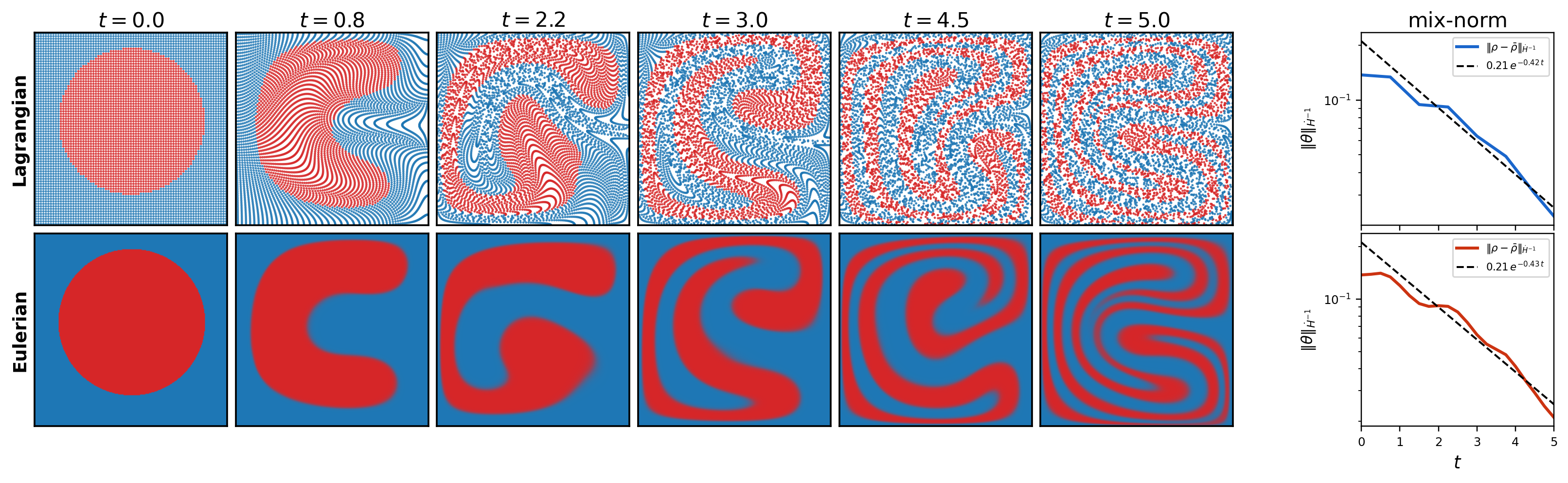}
\caption{Verification for Family 2D-3 (circle).}
\label{fig:eulerian_circle}
\end{figure}

\begin{figure}[!ht]
\centering
\includegraphics[width=\textwidth]{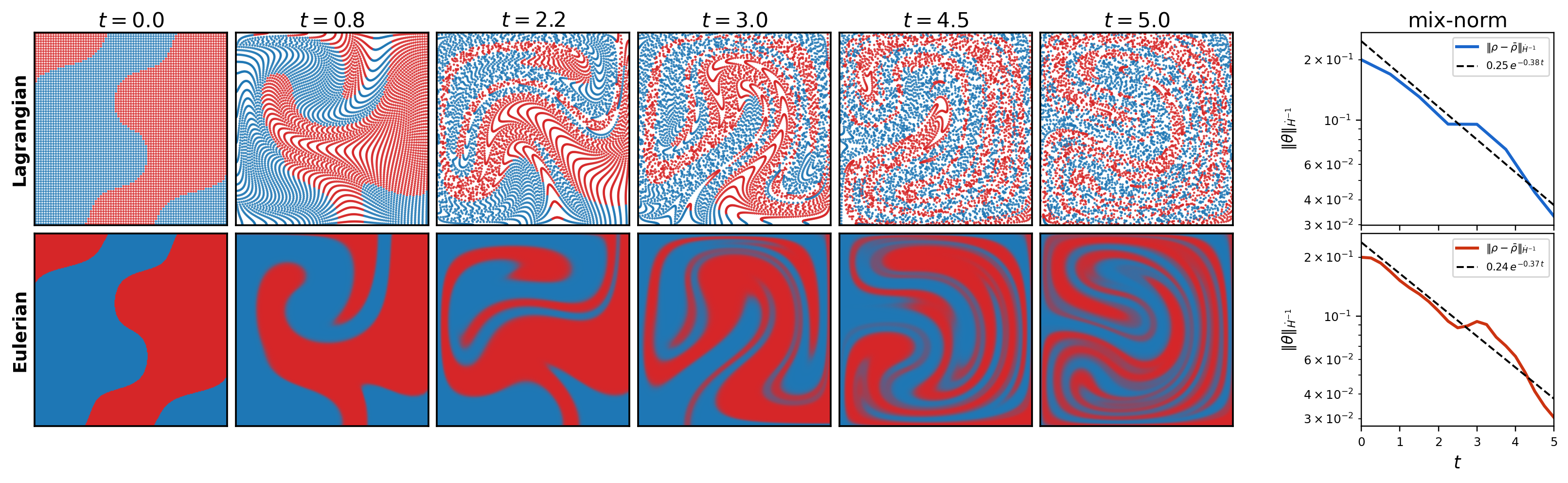}
\caption{Verification for Family 2D-4 (Fourier level set).}
\label{fig:eulerian_level_set}
\end{figure}

On the particle side, the squared mix-norm is evaluated by the spectral estimator \eqref{eq:dye_cost}. Figures~\ref{fig:eulerian_rotation_line}--\ref{fig:eulerian_level_set} compare the particle rollout (top) with the finite-volume rollout (bottom). Both panels use the same initial density, coefficients $\blam^{\mathrm{NN}}(t)$, and time axis. The Eulerian mix-norm is computed by the independent elliptic solve described above.

Figures~\ref{fig:eulerian_rotation_line}--\ref{fig:eulerian_level_set} support two consistency checks. First, the particle snapshots and the gridded scalar snapshots share the same filamentary structure at every recorded time, up to the discretization difference between an empirical particle measure and a piecewise-constant cell average. Second, the exponential decay rates fitted to the two mix-norm curves agree to within the sub-grid filament resolution of the finite-volume solver. The learned velocity field is therefore consistent with the Eulerian transport equation \eqref{eq:transport} on which the problem is originally posed.

%% ====================================================================
\FloatBarrier
\section{Conclusions and perspectives}\label{sec:conclusion}

\textbf{Conclusions.} We proposed an in-context operator-learning framework for optimal mixing of incompressible flows. The framework amortizes the cost of classical PDE-constrained optimization by training a transformer-based operator that reads the two phases of the initial density as point clouds and outputs the time-dependent coefficients of a divergence-free spectral basis. Two design choices make this tractable. First, the Lagrangian reformulation turns the objective into a particle expectation whose gradient with respect to the network parameters is obtained by backpropagation through an explicit RK4 integrator. Second, the spectral basis enforces incompressibility and no-penetration exactly, so no soft penalty on the velocity is required.

The numerical experiments cover four parametric families in 2D and four in 3D, with intrinsic dimensions from $1$ to $9$. In every family the $\Hmone$ mix-norm of the learned flow decays approximately exponentially on the training distribution and generalizes across the held-out parameter hole. Advancing the same velocity field by an independent finite-volume Eulerian solver on a refined grid reproduces the Lagrangian mix-norm decay to within the filament resolution of the solver. Once the operator is trained, computing the optimal control for a new initial density reduces to a single forward pass, in contrast to the per-instance PDE-constrained optimization required by the classical Eulerian approach \cite{hu2019approximating,hu2023feedback}.

\textbf{Perspectives.} Three directions follow. First, multi-family training and architectural extensions aimed at improving cross-family transfer are worth exploring, especially in 3D where the intrinsic complexity of the admissible flows is larger. Second, replacing the uniform spectral basis by an adaptive or learned basis may reduce the required basis size and improve performance on families with strongly anisotropic target geometries. Third, the operator fits naturally into an online feedback loop in which the control is re-predicted from the current density every few time steps. This offers a path toward closed-loop mixing in settings that require real-time decisions, within the Lagrangian framework developed here.

\section*{Acknowledgments}
W. Hu is partially supported by NSF DMS-2111486, DMS-2205117, AFOSR FA9550-23-1-0675 and Humboldt Research Fellowship for Experienced Researchers program from the Alexander von Humboldt Foundation. R. Lai is partially supported by NSF DMS-2401297.

\bibliographystyle{abbrv}
\bibliography{refs_mixing}

\begin{thebibliography}{10}

\bibitem{alberti2016exponential}
G.~Alberti, G.~Crippa, and A.~L. Mazzucato.
\newblock Exponential self-similar mixing by incompressible flows.
\newblock {\em J. Amer. Math. Soc.}, 32(2):445--490, 2019.

\bibitem{balasuriya2015dynamical}
S.~Balasuriya.
\newblock Dynamical systems techniques for enhancing microfluidic mixing.
\newblock {\em J. Micromech. Microeng.}, 25(9):094005, aug 2015.

\bibitem{benner2015survey}
P.~Benner, S.~Gugercin, and K.~Willcox.
\newblock A survey of projection-based model reduction methods for parametric dynamical systems.
\newblock {\em SIAM Rev.}, 57(4):483--531, 2015.

\bibitem{borrvall2003topology}
T.~Borrvall and J.~Petersson.
\newblock Topology optimization of fluids in {S}tokes flow.
\newblock {\em Internat. J. Numer. Methods Fluids}, 41(1):77--107, 2003.

\bibitem{cole2026context}
F.~Cole, D.~Wang, Y.~Chen, Y.~Lu, and R.~Lai.
\newblock In-context operator learning on the space of probability measures.
\newblock {\em arXiv:2601.09979}, 2026.

\bibitem{cortelezzi2008feasibility}
L.~Cortelezzi, A.~Adrover, and M.~Giona.
\newblock Feasibility, efficiency and transportability of short-horizon optimal mixing protocols.
\newblock {\em J. Fluid Mech.}, 597:199--231, 2008.

\bibitem{couchman2010control}
I.~J. Couchman and E.~C. Kerrigan.
\newblock Control of mixing in a stokes’ fluid flow.
\newblock {\em J. Process Control}, 20(10):1103--1115, 2010.

\bibitem{crippa2017cellular}
G.~Crippa and C.~Schulze.
\newblock Cellular mixing with bounded palenstrophy.
\newblock {\em Math. Models Methods Appl. Sci.}, 27(12):2297--2320, 2017.

\bibitem{eggl2020mixing}
M.~F. Eggl and P.~J. Schmid.
\newblock Mixing enhancement in binary fluids using optimised stirring strategies.
\newblock {\em J. Fluid Mech.}, 899:A24, 21, 2020.

\bibitem{eggl2022mixing}
M.~F. Eggl and P.~J. Schmid.
\newblock Mixing by stirring: Optimizing shapes and strategies.
\newblock {\em Phys. Rev. Fluids}, 7:073904, Jul 2022.

\bibitem{feng2026learn}
X.~Feng, L.~Wang, D.~Needell, and R.~Lai.
\newblock Learn to evolve: Self-supervised neural jko operator for wasserstein gradient flow.
\newblock {\em arXiv:2601.05583}, 2026.

\bibitem{franjione1992symmetry}
J.~G. Franjione and J.~M. Ottino.
\newblock Symmetry concepts for the geometric analysis of mixing flows.
\newblock {\em Philosophical Transactions of the Royal Society of London. Series A: Physical and Engineering Sciences}, 338(1650):301--323, 1992.

\bibitem{gao2026self}
S.~Gao, M.~Zhou, and R.~Lai.
\newblock Self-supervised in-context operator learning for stochastic mean-field control.
\newblock {\em arXiv:2608.18282}, 2026.

\bibitem{garcia1989model}
C.~E. García, D.~M. Prett, and M.~Morari.
\newblock Model predictive control: Theory and practice—a survey.
\newblock {\em Automatica}, 25(3):335--348, 1989.

\bibitem{Grebenkov2013Geometrical}
D.~S. Grebenkov and B.-T. Nguyen.
\newblock Geometrical structure of {L}aplacian eigenfunctions.
\newblock {\em SIAM Rev.}, 55(4):601--667, 2013.

\bibitem{gubanov2010towards}
O.~Gubanov and L.~Cortelezzi.
\newblock Towards the design of an optimal mixer.
\newblock {\em J. Fluid Mech.}, 651:27--53, 2010.

\bibitem{hu2025cellular}
W.~Hu, M.-J. Lai, and H.-N. Wu.
\newblock Cellular flow control design for mixing based on the least action principle.
\newblock {\em arXiv:2510.22703}, 2025.

\bibitem{hu2026structure}
W.~Hu, Z.~Li, Y.~Zhang, and E.~Zuazua.
\newblock A structure-preserving numerical scheme for optimal control and design of mixing in incompressible flows.
\newblock {\em arXiv:2601.06294}, 2026.

\bibitem{hu2023feedback}
W.~Hu, C.~N. Rautenberg, and X.~Zheng.
\newblock Feedback control for fluid mixing via advection.
\newblock {\em J. Differential Equations}, 374:126--153, 2023.

\bibitem{hu2019approximating}
W.~Hu and J.~Wu.
\newblock An approximating approach for boundary control of optimal mixing via {N}avier-{S}tokes flows.
\newblock {\em J. Differential Equations}, 267(10):5809--5850, 2019.

\bibitem{huang2025unsupervised}
H.~Huang and R.~Lai.
\newblock Unsupervised solution operator learning for mean-field games.
\newblock {\em J. Comput. Phys.}, 537:Paper No. 114057, 22, 2025.

\bibitem{kovachki2023neural}
N.~Kovachki, Z.~Li, B.~Liu, K.~Azizzadenesheli, K.~Bhattacharya, A.~Stuart, and A.~Anandkumar.
\newblock Neural operator: learning maps between function spaces with applications to {PDE}s.
\newblock {\em J. Mach. Learn. Res.}, 24:Paper No. [89], 97, 2023.

\bibitem{LeVeque2002FVM}
R.~J. LeVeque.
\newblock {\em Finite volume methods for hyperbolic problems}.
\newblock Cambridge Texts in Applied Mathematics. Cambridge University Press, Cambridge, 2002.

\bibitem{li2020fourier}
Z.~Li, N.~Kovachki, K.~Azizzadenesheli, B.~Liu, K.~Bhattacharya, A.~Stuart, and A.~Anandkumar.
\newblock Fourier neural operator for parametric partial differential equations.
\newblock In {\em International Conference on Learning Representations}, 2021.

\bibitem{li2026hamiltonian}
Z.~Li and E.~Zuazua.
\newblock Hamiltonian interface dynamics for reduced-order optimization of incompressible mixing.
\newblock {\em arXiv:2605.04688}, 2026.

\bibitem{lin2011optimal}
Z.~Lin, J.-L. Thiffeault, and C.~R. Doering.
\newblock Optimal stirring strategies for passive scalar mixing.
\newblock {\em J. Fluid Mech.}, 675:465--476, 2011.

\bibitem{liu1994quantification}
M.~Liu, F.~Muzzio, and R.~Peskin.
\newblock Quantification of mixing in aperiodic chaotic flows.
\newblock {\em Chaos Solitons Fractals}, 4(6):869--893, 1994.
\newblock Special Issue: Chaos Applied to Fluid Mixing.

\bibitem{lu2021deeponet}
L.~Lu, P.~Jin, G.~Pang, Z.~Zhang, and G.~E. Karniadakis.
\newblock Learning nonlinear operators via {DeepONet} based on the universal approximation theorem of operators.
\newblock {\em Nat. Mach. Intell.}, 3(3):218--229, 2021.

\bibitem{lunasin2012optimal}
E.~Lunasin, Z.~Lin, A.~Novikov, A.~Mazzucato, and C.~R. Doering.
\newblock Optimal mixing and optimal stirring for fixed energy, fixed power, or fixed palenstrophy flows.
\newblock {\em J. Math. Phys.}, 53(11):115611, 15, 2012.

\bibitem{mathew2007optimal}
G.~Mathew, I.~Mezi\'c, S.~Grivopoulos, U.~Vaidya, and L.~Petzold.
\newblock Optimal control of mixing in {S}tokes fluid flows.
\newblock {\em J. Fluid Mech.}, 580:261--281, 2007.

\bibitem{mathew2005multiscale}
G.~Mathew, I.~Mezi\'c, and L.~Petzold.
\newblock A multiscale measure for mixing.
\newblock {\em Phys. D}, 211(1-2):23--46, 2005.

\bibitem{sharma1997control}
A.~Sharma and N.~Gupte.
\newblock Control methods for problems of mixing and coherence in chaotic maps and flows.
\newblock {\em Pramana}, 48(1):231--248, Jan 1997.

\bibitem{smith2013uncertainty}
R.~C. Smith.
\newblock {\em Uncertainty quantification}, volume~12 of {\em Computational Science \& Engineering}.
\newblock Society for Industrial and Applied Mathematics (SIAM), Philadelphia, PA, 2014.
\newblock Theory, implementation, and applications.

\bibitem{stroock2002chaotic}
A.~D. Stroock, S.~K.~W. Dertinger, A.~Ajdari, I.~Mezic, H.~A. Stone, and G.~M. Whitesides.
\newblock Chaotic mixer for microchannels.
\newblock {\em Science}, 295(5555):647--651, 2002.

\bibitem{thiffeault2012using}
J.-L. Thiffeault.
\newblock Using multiscale norms to quantify mixing and transport.
\newblock {\em Nonlinearity}, 25(2):R1--R44, 2012.

\bibitem{vaswani2017attention}
A.~Vaswani, N.~Shazeer, N.~Parmar, J.~Uszkoreit, L.~Jones, A.~N. Gomez, {\L}.~Kaiser, and I.~Polosukhin.
\newblock Attention is all you need.
\newblock In {\em Advances in Neural Information Processing Systems (NeurIPS)}, volume~30, pages 5998--6008, 2017.

\bibitem{vikhansky2002control}
A.~Vikhansky.
\newblock Control of stretching rate in time-periodic chaotic flows.
\newblock {\em Phys. Fluids}, 14(8):2752--2756, 2002.

\bibitem{vikhansky2002enhancement}
A.~Vikhansky.
\newblock Enhancement of laminar mixing by optimal control methods.
\newblock {\em Chem. Eng. Sci.}, 57(14):2719--2725, 2002.

\bibitem{villermaux2019mixing}
E.~Villermaux.
\newblock Mixing versus stirring.
\newblock In {\em Annual review of fluid mechanics. {V}ol. 51}, volume~51 of {\em Annu. Rev. Fluid Mech.}, pages 245--273. Annual Reviews, Palo Alto, CA, 2019.

\bibitem{incontext2023}
L.~Yang, S.~Liu, T.~Meng, and S.~J. Osher.
\newblock In-context operator learning with data prompts for differential equation problems.
\newblock {\em Proc. Natl. Acad. Sci. USA}, 120(39):Paper No. e2310142120, 10, 2023.

\bibitem{zaheer2017deepsets}
M.~Zaheer, S.~Kottur, S.~Ravanbakhsh, B.~Poczos, R.~R. Salakhutdinov, and A.~J. Smola.
\newblock Deep sets.
\newblock {\em Advances in Neural Information Processing Systems}, 30, 2017.

\bibitem{zheng2023numerical}
X.~Zheng, W.~Hu, and J.~Wu.
\newblock Numerical algorithms and simulations of boundary dynamic control for optimal mixing in unsteady {S}tokes flows.
\newblock {\em Comput. Methods Appl. Mech. Engrg.}, 417:Paper No. 116455, 24, 2023.

\end{thebibliography}

\end{document}